\documentclass[10pt,openbib]{article}
\usepackage[T1]{fontenc}
\usepackage{amsmath}
\usepackage{amsfonts}
\usepackage{amssymb}
\usepackage{graphicx}
\newcommand{\N}{\mathbb{N}}
\newcommand{\R}{\mathbb{R}} 

\newcommand{\Q}{\mathbb{Q}}

\newcommand{\Z}{\mathbb{Z}}

\newcommand{\dee}{\mathop{\! \,\mathrm{d} \!}\nolimits}
\newcommand{\comp}{\raisebox{0pt}{$\scriptstyle\circ \, $}}
\newcommand{\setrule}{\, \rule[-4pt]{.5pt}{13pt}\, }

\newcommand{\dbydt}{\mbox{${\displaystyle \frac{\dee}{\dee t}}
\rule[-10pt]{.5pt}{25pt} \raisebox{-10pt}{$\, {\scriptstyle t=0}$}$}}

\newcommand{\vvee}{\mbox{\tiny $\vee $}}

\begin{document}
\begin{center}
{\Large \textbf{Differentiable spaces that are subcartesian}} \\[0pt]
\mbox{}\vspace{-0.2in} \\[0pt]
\mbox{}\\[0pt]
Richard Cushman and J\k{e}drzej \'{S}niatycki\footnotemark 
\end{center}

\footnotetext{version: \today} 
\footnotetext{
Department of Mathematics and Statistics, University of Calgary \\ 
e-mail: rcushman@ucalgary.ca  and sniatycki@ucalgary.ca \\
}
\bigskip 

\begin{description}
\item[Abstract] We show that the differential structure of the orbit space 
of a proper action of a Lie group on an a smooth manifold is continuously reflexive. This 
implies that the orbit space is a differentiable space in the sense of Smith, which ensures 
that the orbit space has an exterior algebra
of differential forms, which satisfies Smith's version of de Rham's theorem. Because 
the orbit space is a locally closed subcartesian space, it has vector fields and their flows.
\end{description}

\section{Introduction}


The aim of this paper is to show that the differential structure of the orbit space 
$M/G$ of a proper action of a Lie group $G$ on an a smooth manifold $M$ is 
continuously reflexive. This implies that $M/G$ a differentiable space in the sense of Smith. 
This ensures that $M/G$ has an exterior algebra of
differential forms, which satisfies Smith's version of de Rham's theorem. Because $M/G$ 
is a locally closed subcartesian space, it has vector fields and their flows. \smallskip

The concepts of a differentiable structure and a differentiable space were 
introduced by Smith in 1966 \cite{smith}. They allowed him to define the notions of
a differential form and exterior derivative, which he used to prove a
version of de Rham's theorem. Smith's paper has not received the attention
of mathematical community that it deserves. In section 5 we give a
comprehensive review of results in Smith \cite{smith}. \smallskip 

Subcartesian spaces and differential spaces were introduced in 1967 by
Aronszajn \cite{aronszajn} and Sikorski \cite{sikorski67}, respectively.
Subcartesian spaces are Hausdorff spaces locally diffeomorphic to arbitrary
subsets of ${\mathbb{R}^{n}}$, where their differential structure is given
by an atlas of singular charts. The differential structure of a differential
space is given by its ring of smooth functions. In 1973, Walczak \cite{walczak} showed that subcartesian spaces can be treated as a subcategory of
differential spaces. This result allows for a discussion of the geometry of
a subcartesian space $S$ in terms of its differential structure 
$C^{\infty}(S)$. In particular, derivations of $C^{\infty }(S)$ give rise to the
notion of vector fields on $S$ and their flows. We summarize the aspects of
the theory of subcartesian spaces needed here in section 2 following \cite{sniatycki} supplemented by \cite{cushman-sniatycki19}. \smallskip 

We introduce the notion of a continuously reflexive differential structure,
which is modeled on the notion of a reflexive differential structure
discussed in Batubenge et al. \cite{batubenge-iglesias-karshon-watts}. We show that, if the differential structure $C^{\infty }(S)$ of a subcartesian space $S$ is continuously reflexive, then
it is a differentiable structure on $S$ in the sense of Smith. This
ensures that a subcartesian space $S$ with continuously reflexive
differential structure $C^{\infty }(S)$ has an exterior algebra of
differential forms, which satisfies Smith's version of de Rham's theorem. \smallskip

The category of locally closed subcartesian differential spaces whose 
differential structure is continuously reflexive and smooth maps between them 
contains the category of orbit spaces of a proper action on a smooth manifold and 
smooth maps between them. In more detail, let $\Phi $ be a proper action of a Lie group $G$ on a connected smooth manifold $M$. By a result of Palais \cite{palais61} the $G$-action $\Phi $ has a slice $S_m$ at each point $m \in M$. We denote the space of smooth $G$-invariant functions on $M$ by ${C^{\infty}(M)}^G$. Let  
$M/G = \{ G \cdot m \setrule \, m \in M \}$ be the space of $G$-orbits in $M$ with $G$-orbit map 
$\pi : M \rightarrow M/G: m \mapsto G\cdot m$. Classically it is known 
that $M/G$ with the quotient topology is a Hausdorff topological manifold, 
see \cite[chpt VII]{cushman-bates15}. The space $C^{\infty}(M/G)= 
\{ f : \overline{M} \rightarrow \R \setrule {\pi }^{\ast }f \in {C^{\infty}(M)}^G \} $ is a 
differential structure on $M/G$, see Cushman and \'{S}niatycki \cite{cushman-sniatycki}. 
The differential space topology on $M/G$ is equivalent to the quotient topology, see 
\cite[prop. 4.3.1, p.68]{sniatycki}. For each $m \in M$ let $G_m = \{ g \in G \setrule \, {\Phi}_g(m) 
= m \}$. Since the differential spaces $(S_m/G_m, (C^{\infty}(T_mS_m)^{G_m})$ and $(M/G, C^{\infty}(M)^G)$ are locally diffeomorphic, see \cite{cushman-sniatycki}, it follows that the differential space $(M/G,C^{\infty}(M/G))$ is subcartesian and is locally closed, 
see \cite[lem. 4.6]{duistermaat}. Hence $M/G$ is paracompact. In section six we show that the differential structure $C^{\infty}(M/G)$ is continuously reflexive. Thus $M/G$ is a differentiable space in the sense of Smith.

\section{Subcartesian differential spaces}

A \emph{differential structure} on a topological space 
$S$ is a family $C^{\infty}(S)$ of real valued functions on $S$ satisfying the following conditions: \medskip 

1. \parbox[t]{4in}{The collection $\{ f^{-1}(I) \setrule \, f \in C^{\infty}(S) \, 
\mbox{and $I$ is open interval in $\R $} \} $  is a subbasis for the topology of $S$.} \smallskip 

2. \parbox[t]{4in}{For every $n \in \N$, if $f_1, \ldots , f_n \in C^{\infty}(S)$ and 
$F \in C^{\infty}({\R }^n)$, then $F(f_1,\ldots ,f_n) \in C^{\infty}(S)$.} \smallskip 

3. \parbox[t]{4in}{If $f :S \rightarrow \R$ is a function such that for every $x \in S$ there exists an open neighbourhood $U$ of $x$ and a function $f_x \in C^{\infty}(S)$ satisfying 
$f_x|_U = f|_U$, then $f \in C^{\infty} (S)$. Here, the vertical bar denotes restriction.} \bigskip 

\noindent A function $f \in C^{\infty}(S)$ is called a \emph{smooth} function on $S$. It follows from 
condition 1, that a smooth function on $S$ is continuous. Condition 2 with $F(f_1, f_2) = 
a f_1 +b f_2$, where $a$, $b \in \R $ shows that $C^{\infty}(S)$ is a real vector space. Taking 
$F(f_1, f_2) = f_1\cdot f_2$ shows that $(C^{\infty}(S), \cdot )$ is a commutative associative 
algebra with unit under multiplication $\cdot $. \medskip 

In his original definition Sikorski \cite{sikorski67} defined $C^{\infty}(S)$ to be a family of functions 
on $S$ satisfying condition 2. He then used condition 1 to define a topology on $S$. Finally, he 
imposed condition 3 as a consistency condition. \medskip  

A topological space $S$ endowed with a differential structure $C^{\infty}(S)$ is called a 
\emph{differential space}. \medskip 

A map $\varphi : R \rightarrow S$ between diferential spaces is 
\emph{smooth} if ${\varphi }^{\ast }f = f \comp \varphi \in C^{\infty}(R)$ for every $f \in C^{\infty}(S)$. 
The map $\varphi $ is a \emph{diffeomorphism} if it is invertible and 
${\varphi }^{-1}$ is smooth. \medskip 

\noindent \textbf{Proposition 2.1} A smooth map between differential spaces is continuous. \medskip 

\noindent \textbf{Proof.} See proposition 2.1.5 in \cite{sniatycki}. \hfill $\square $ \medskip

An alternative way of constructing a differential structure on a \emph{set} $S$, used by Sikorski 
\cite{sikorski72}, goes as follows. Let $\mathcal{F}$ be a family of real valued functions on 
$S$. \linebreak 
Endow $S$ with the topology generated by the subbasis $\{ f^{-1}(I) \, \mbox{$f \in \mathcal{F}$ 
and $I$ an} \linebreak 
\mbox{open interval in $\R $} \} $. Define $C^{\infty}(S)$ by saying that $h \in C^{\infty}(S)$ if and only if for each $x \in S$ there is an  
open subset $U$ of $S$, functions $f_1, \ldots , f_n \in \mathcal{F}$, and 
$F \in C^{\infty}({\R }^n)$ such that $h|_U = F(f_1, \ldots , f_n)|_U$. Clearly $\mathcal{F} \subseteq 
C^{\infty}(S)$. Moreover, \medskip 

\vspace{-.15in}\noindent \textbf{Theorem 2.2} The family $C^{\infty}(S)$ is a differential structure on $S$. \medskip 

\noindent \textbf{Proof.} See theorem 2.1.7 in \cite{sniatycki}. \hfill $\square $ \medskip 

We refer to $C^{\infty}(S)$ as the differential structure on $S$ \emph{generated} by 
the family $C^{\infty}(S)$. \medskip 

\noindent \textbf{Proposition 2.3} Let $M$ be a smooth manifold. Then $C^{\infty}(M)$ is a 
differential structure. \medskip 

\noindent \textbf{Proof.} To see this we check that the defining properties of a differential 
structure hold. First, since every function $f \in C^{\infty}(M)$ is continuous, the set 
$\{ f^{-1}(I) \subseteq M \setrule \, f \in C^{\infty}(M) \, \, \mbox{and $I$ an open interval in $\R $} \} $ 
is a subbasis for the manifold topology of $M$. Second, for every $n \in \N$ with 
$f_1, \ldots, f_n \in C^{\infty}(M)$ and $F \in C^{\infty}({\R }^n)$ we see that 
$F(f_1, \ldots , f_n) \in C^{\infty}(M)$. Third, suppose that $f$ is a function on $M$ 
such that for every $m \in M$ there is an open neighborhood $U$ of $m$ in $M$ and 
a function $f_m \in C^{\infty}(M)$ such that $f|_U = f_m|_U$. Choosing $U$ to be the domain 
of a chart of $M$ at $m$, it follows that $f$ is a smooth function on $M$. Consequently, 
$C^{\infty}(M)$ is a differential structure on $M$. \hfill $\square $ \medskip 

Let $R$ be a differential space with differential structure $C^{\infty}(R)$. Let $S$ be a 
\emph{subset} of $R$ endowed with the subspace topology, in which open subsets of $S$ 
are of the form $S \cap U$, where $U$ is an open subset of $R$. Let 
$C^{\infty}(R)$ be the differential structure generated by the restriction to $R$ of smooth functions on $S$. 
We say that the differential space $(S, C^{\infty}(S))$ is a \emph{differential subspace} 
of the differential space $(R, C^{\infty}(R))$. \medskip 

\noindent \textbf{Proposition 2.4} The space $C^{\infty}(R)|_S$ of restrictions to $S \subseteq R$ 
of smooth functions on $R$ generates a differential structure $C^{\infty}(S)$ on $S$ such that 
the differential space topology on $S$ coincides with its subspace topology. The inclusion map 
$\iota : S \rightarrow R$ is smooth. \medskip 

\noindent \textbf{Proof.} See proposition 2.1.8 in \cite{sniatycki}. \hfill $\square $ \medskip 

We can characterize a smooth $n$-dimensional Hausdorff manifold as a Hausdorff differential 
space $S$ such that every point $x \in S$ has a neighborhood diffeomorphic to an open subset 
$V$ of ${\R }^n$. Here the differential structures on $U$ and $V$ are generated by the restriction of smooth functions on $S$ to $U$ and on ${\R }^n$ to $V$, respectively. \medskip 

We weaken the above definition of manifold as follows.  A differential space $(S, C^{\infty}(S))$ is \emph{subcartesian} if its topology is Hausdorff and for every point $x \in S$ there is an open neighborhood $U$, which is diffeomorphic to a \emph{subset} $V$ of some ${\R }^n$. Observe that $V$ is an arbitrary subset of ${\R }^n$ and that $n$ may depend on $x \in S$. More precisely, a differential space $(S, C^{\infty}(S))$ is subcartesian if and only 
if \medskip 

\par \noindent \begin{tabular}{l}
1. \parbox[t]{4.2in}{the topology on $S$ defined by its differential structure 
$C^{\infty}(S)$ is Hausdorff;} \\
\end{tabular}
\par \noindent 
\begin{tabular}{l}
2. \parbox[t]{4.2in}{for every $x \in S$ there is an open subset $U_x$ of $S$ 
containing $x$ and a diffeomorphism 
${\varphi }_x: \big( U_x, C^{\infty}(U_x) \big) \rightarrow 
\big( V_x = {\varphi }_x(U_x), C^{\infty}(V_x) \big) $ 
of differential spaces, where $V_x \subseteq {\R }^n$, ${\iota }_{V_x}: V_x \rightarrow {\R }^n$ is the inclusion map and $C^{\infty}(V_x) = {\iota }^{\ast }_{V_x}(C^{\infty}({\R }^n))$.}\footnotemark
\end{tabular}  
\footnotetext{In other words, $f \in C^{\infty}(V_x)$ if and only if there is an 
$F \in C^{\infty}({\R }^n)$ such that $F|_{V_x} = f$.}  \medskip 

\noindent A diffeomorphism of an open subset of $S$ onto a subset of some ${\R }^n$ is  
called a $\emph{chart}$ on $S$. The family of all charts is the \emph{complete atlas} on 
$S$. An example of a subcartesian differential space is an arbitrary subset $W$ of Euclidean space 
${\R}^n$ with subset topology and differential structure $C^{\infty}({\R }^n)|_W$ consisting of restrictions of smooth functions on ${\R}^n$ to $W$. We refer to a differential space, which is subcartesian, 
as a \emph{subcartesian differential space}, because subcartesian spaces with charts is an active 
area of research.

\section{Vector fields on a subcartesian space}

\subsection{Tangent bundle}

Let $(S, C^{\infty}(S))$ be a differential space. A \emph{derivation} of $C^{\infty}(S)$ 
\emph{at} $x \in S$ is a real linear mapping $v_x: C^{\infty}(S) \rightarrow \R $ 
such that 
\begin{equation}
v_x(f_1 f_2) = v_x(f_1)f_2(x) + f_1(x) v_x(f_2) \quad \mbox{for every $f_1$, $f_2 \in C^{\infty}(S)$}. 
\label{eq-nsec4ss1one}
\end{equation}
Let $T_xS$, the \emph{tangent space} to $S$ \emph{at} $x$, be the set of all derivations 
of $C^{\infty}(S)$ at $x$. Let $TS =\{ v = (x, v_x) \setrule \, x\in S \, \, \& \, \, v_x \in T_xS \} $ be the \emph{tangent bundle} 
of $S$ with bundle \emph{projection map} ${\tau }_S: TS \rightarrow S: v = (x, v_x) \mapsto x$, 
which assigns to every derivation $v_x \in T_xS$ at $x$ the point $x$. The differential 
structure $C^{\infty}(TS)$ on $TS$ is generated by the family of functions 
$\{ {\tau }^{\ast }_Sf, \dee f \setrule \, f \in C^{\infty}(S) \}$, where 
$\dee f: TS \rightarrow \R : v = (x, v_x) \mapsto v_x(f)$. The tangent bundle projection map 
${\tau }_S : (TS, C^{\infty}(TS)) \rightarrow (S, C^{\infty}(S))$ is a smooth map of differential spaces. \medskip 

\noindent \textbf{Lemma 3.1.1} If $\varphi : S \rightarrow R$ is a smooth map of differential 
spaces, then 
\begin{displaymath}
T\varphi :TS \rightarrow TR: v = (x, v_x) \mapsto \big( \varphi (x), T_x{\varphi }(v_x) \big) , 
\end{displaymath}
where $\big(T_x{\varphi }(v_x) \big)f = v_x({\varphi }^{\ast }f)$ for every $f \in C^{\infty}(R)$, is a smooth map of differential spaces. \medskip 

\noindent \textbf{Proof.} It suffices to show that  
$(T\varphi )^{\ast }{\mathcal{F}}_R \subseteq {\mathcal{F}}_S$, 
where ${\mathcal{F}}_R = \{ {\tau }^{\ast }_R f, \dee f \setrule \, f \in C^{\infty}(R) \}$ and 
${\mathcal{F}}_S = \{ {\tau }^{\ast }_Sg, \dee g \setrule \, g \in C^{\infty}(S) \}$. 
Since $T\varphi $ sends $T_xS$ to $T_{\varphi (x)}R$ for every $x \in S$, we have 
${\tau }_R \comp T\varphi = \varphi \comp {\tau }_S$. So for every $f \in C^{\infty}(R)$ we get 
\begin{displaymath}
(T\varphi )^{\ast }({\tau }^{\ast }_Rf) = ({\tau }_R \comp T{\varphi })^{\ast }f 
= (\varphi \comp {\tau }_S)^{\ast}f = {\tau }^{\ast }_S({\varphi }^{\ast }f) = 
{\tau }^{\ast }_S(g), 
\end{displaymath}
for $g = {\varphi }^{\ast}f \in C^{\infty}(S)$. And for every $v = (x,v_x) \in TS$ we have 
\begin{align}
(T\varphi )^{\ast }(\dee f)(v) &  = \big( T_x\varphi (v_x)(f) \big) (x) 
= \big( v_{\varphi (x)}({\varphi }^{\ast }f)\big) (x) = \dee \, ({\varphi }^{\ast }f)(v)  , 
\notag 
\end{align}
that is, $(T\varphi )^{\ast} (\dee f) = \dee \, ({\varphi }^{\ast }f) = \dee g$, for 
$g = {\varphi }^{\ast } f \in C^{\infty}(S)$. Consequently, $(T\varphi )^{\ast }{\mathcal{F}}_R 
\subseteq {\mathcal{F}}_S$, as desired. \hfill $\square $ \medskip

\noindent \textbf{Corollary 3.1.1A} Suppose that $\varphi : R \rightarrow S$ and $\psi : S \rightarrow 
U$ are smooth maps of differential spaces, then $T(\psi \comp \varphi ) =T\psi \comp T\varphi $. 
\medskip 

\noindent \textbf{Proof.} Let $v_x \in T_xR$. Then $v_{\varphi (x)} = T_x{\varphi }(v_x) 
\in T_{\varphi (x)}S$. Also $v_{\psi (\varphi (x))} = 
T_{\varphi (x)}\psi (v_{\varphi (x)}) \in T_{\psi (\varphi (x))} U$. 
So $T_{\psi (\varphi (x))}\psi \comp T_{\varphi (x)}\varphi $ sends $v_x$ to 
$v_{(\psi \comp \varphi )(x)}$, that is, 
$T_{\psi (\varphi (x))}\psi \comp T_{\varphi (x)}\varphi = T_x(\psi \comp \varphi )$ for 
every $x \in R$. Thus $T (\psi \comp \varphi ) = T\psi \comp T \varphi $. \hfill $\square $ \medskip  

\noindent \textbf{Corollary 3.1.1B} If $\varphi : S \rightarrow R$ is a diffeomorphism of 
differential spaces, then $T\varphi :TS \rightarrow TR$ is a diffeomorphism of 
differential spaces. \medskip 

\noindent \textbf{Proof.} This follows immediately because $(T\varphi )^{-1} = 
T{\varphi }^{-1}$. \hfill $\square $ \medskip 

\noindent \textbf{Corollary 3.1.1C} If $(R, C^{\infty}(R))$ is a differential subspace of 
$(S, C^{\infty}(S))$, then $(TR, C^{\infty}(R))$ is a differential subspace of 
$(TS, C^{\infty}(TS))$. In particular, $TR \subseteq TS$. \medskip 

\noindent \textbf{Proof.} Because $R$ is a differential subspace of $S$, the inclusion 
map ${\iota }_R: R \rightarrow S$ is a smooth map. Hence its tangent 
\begin{displaymath}
T{\iota }_R: TR \rightarrow TS: (x, v_x) \mapsto (y, v_y) = 
({\iota }_R(x), (T_x{\iota }_R) v_x) 
\end{displaymath} 
is a smooth map. $T{\iota }_R$ is the inclusion map 
because $T_x{\iota }_R: T_x R \rightarrow T_{{\iota }_R (x)}S$ 
is the inclusion map for every $x \in R$. Hence $(TR, C^{\infty}(TR))$ is a differential 
subspace of $(TS, C^{\infty}(TS))$. In particular, $TR \subseteq TS$. \hfill $\square $ \medskip 

\noindent \textbf{Proposition 3.1.2} If $(S, C^{\infty}(S))$ is a subcartesian differential 
space, then its tangent bundle $(TS, C^{\infty}(TS))$ is a 
subcartesian differential space.  \medskip 

\noindent \textbf{Proof.} It suffices to show that the smooth functions 
in $\mathcal{F} = \{ {\tau }^{\ast}_Sf, \, \dee f \setrule \, f \in C^{\infty}(S) \} $ separate 
points of $TS$. Let $v = (x,v_x)$ and $w = (y, w_y)$ be two distinct points in 
$TS$. Suppose that $x = {\tau }_S(v) \ne {\tau }_S(w) = y$. Since $S$ is subcartesian, 
the topology on $S$ determined by $C^{\infty}(S)$ is Hausdorff. So there are 
disjoint open neighborhoods $V$ of $x$ and $W$ of $y$ in $S$. Since the 
bundle projection map ${\tau }_S:TS \rightarrow S$ is smooth, it is continuous. 
Hence ${\tau }^{-1}_S(V)$ and ${\tau }^{-1}_S(W)$ are open neighborhoods in $TS$ of 
$v$ and $w$, respectively. Also ${\tau }^{-1}_S(V) \cap {\tau }^{-1}_S(W) = \varnothing $, 
for if $u \in {\tau }^{-1}_S(V) \cap {\tau }^{-1}_S(W)$, then ${\tau }_S(u) \in V \cap W$, 
which is a contradiction. Now suppose that $x =y$. Since $v \ne w$, it follows that 
$v_x \ne w_x$. Thus for some $f \in C^{\infty}(S)$ we have 
\begin{displaymath}
a = \dee f(v) = \dee f(x)v_x = v_x(f) \ne w_x(f) = \dee f(x)w_x = \dee f(w) = b, 
\end{displaymath}
where $\dee f(x):T_xS \rightarrow \R :v_x \mapsto v_x(f)$. 
Since either $a$ or $b$ is nonzero, the linear map $\dee f(x)$ is nonzero and hence is surjective. Consequently, the map $\dee f: TS \rightarrow \R $ is surjective. Suppose that $a < b$. The case when $b < a$ is handled similarly and is omitted. 
We can choose $c < d$ so that $a < c < d < b$. For example let $c = a +\frac{1}{4}(b-a)$ and 
$d = b - \frac{1}{4}(b-a)$. Let $V = (\dee f)^{-1}(-\infty, c)$ and $W = (\dee f)^{-1}(d, \infty )$. 
Then $V$ and $W$ are open subsets of $TS$, since $\dee f$ is a smooth, and hence a continuous, function on $TS$. Moreover, $V \cap W = \varnothing $, because 
$(-\infty ,c) \cap (d, \infty) = \varnothing $. Also $v \in V$ and $w \in W$, since 
$\dee f(v) = \dee f(x)v_x = a \in (-\infty ,c)$ and $\dee f(w) = \dee f(x)w_x = b \in (d , \infty)$. 
Thus the topology on $TS$ generated by the functions in $\mathcal{F}$, and hence in 
$C^{\infty}(TS)$, is Hausdorff. \hfill $\square $ \medskip  

A \emph{section} of the tangent bundle projection map ${\tau }_S: TS \rightarrow S$ is a 
smooth map $\xi : S \rightarrow TS$ of differential spaces such that ${\tau }_S \comp \xi 
= {\mathrm{id}}_S$. A real linear map $X: C^{\infty}(S) \rightarrow C^{\infty}(S)$ is 
a \emph{derivation} on the differential space $(S, C^{\infty}(S))$ if it satisfies Leibniz's rule
\begin{equation}
X(f_1f_2) = X(f_1)f_2 + f_1X(f_2) \quad \mbox{for every $f_1$, $f_2 \in C^{\infty}(S)$}. 
\label{eq-ns4ss1two}
\end{equation}

\noindent \textbf{Proposition 3.1.3} Every derivation $X$ of $C^{\infty}(S)$ defines 
a section $X: S \rightarrow TS: x \mapsto X(x)$ of ${\tau }_S$, where $X(x)f = (Xf)(x)$ for every 
$f \in C^{\infty}(S)$. \medskip 

\noindent \textbf{Proof.} See \cite[prop. 6]{cushman-sniatycki19}. \hfill $\square $ 

\subsection{Integral curves of derivations}

Let $c: I \subseteq \R \rightarrow S$ be a smooth map of an interval $I$ containing $0$ into the 
differential space $S$. In other words, $c$ is a smooth map of the differential space 
$(I, C^{\infty}(I))$ into the differential space $(S, C^{\infty}(S))$ with $0 \in I$. 
We say that $c$ is an \emph{integral curve} of the \emph{derivation} 
$X$ of $C^{\infty}(S)$ \emph{starting at} $x_0 \in S$ if $c(0) = x_0$ and 
\begin{displaymath}
\frac{\dee }{\dee t}f\big( c(t) \big) = X(f)\big( c(t) \big) \quad 
\mbox{for every $f \in C^{\infty}(S)$ and every $t \in I$} . 
\end{displaymath}
So $c: I \subseteq \R \rightarrow S$ is an integral curve of $X$ if 
$(Tc)(t) = (X\comp c)(t)$. The following example shows that it is convenient to allow the 
interval $I$ to be a single point. \medskip 

\noindent \textbf{Example 3.2.1} Let $\Q $ be the set of rational numbers in $\R $. 
Since $\Q $ is dense in $\R $, every 
derivation of $C^{\infty}(\R )$ induces a derivation of $C^{\infty}(\Q)$. Let 
$X$ be the derivation of $C^{\infty}(\Q)$ induced by the derivation $\frac{\dee }{\dee x}$ of 
$C^{\infty}(\R )$. In other words, for every $f \in C^{\infty}(\Q)$ and every $x_0 \in \Q $ 
\begin{displaymath}
(Xf)(x_0) = \lim_{x \rightarrow x_0} \frac{f(x) - f(x_0)}{x-x_0}, 
\end{displaymath}
where the limit is taken over $x \in \Q $. Since no two distinct points in $Q$ can be 
joined by a continuous curve in $\Q $, it follows that for each $x_0 \in \Q $ the 
tangent vector $X(x_0) \in T_{x_0}\Q $ is the maximal integral curve of $X$ through $x_0$. 
\hfill $\square $ \medskip 

Let $(S, C^{\infty}(S))$ be a differential space. Let $X$ be a derivation of 
$C^{\infty}(S)$. A \emph{lifted integral curve} of $X$ \emph{starting at} $x_0 \in S$ is a 
map $\gamma : I \subseteq \R \rightarrow TS$ such that $\gamma (0) = X(x_0)$ and 
\begin{displaymath}
\frac{\dee }{\dee t} f \big( {\tau}_S (\gamma (t)) \big) = X(f)\big( {\tau}_S (\gamma (t)) \big) , 
\end{displaymath}
for every $f \in C^{\infty}(S)$ and every $t \in I$ provided that $I \ne \{ 0 \}$. \medskip 

\noindent The following example shows that quite ordinary differential spaces 
can have a derivation with an integral curve whose domain is a single point. \medskip 

\noindent \textbf{Example 3.2.2} Consider the function 
$h(x) =${\tiny $\left\{ \begin{array}{cl}\hspace{-4pt} {\mathrm{e}}^{(-x^{-2})}, & \hspace{-8pt} \mbox{if $x \ne 0$} \\
\hspace{-7pt} 0, &\hspace{-8pt} \mbox{if $x=0$.} \end{array} \right. $} Then $h$ is smooth and it and all of its derivatives vanish at $x=0$. Let 
\begin{displaymath}
S = \{ (x,y)\in {\R }^2 \setrule \, H(x,y) = y^2-h(x)y = 0 \} .
\end{displaymath} 
The point $(0,0)$ is an isolated singular point of $S$, since 
$\dee H(x,y) = (-h^{\prime}(x)y, 2y - h(x))$ 
is equal to $(0,0)$ only when $(x,y) = (0,0)$.  \medskip 

\noindent Consider the map $X = \frac{\partial }{\partial x} + v(x) \frac{\partial }{\partial y}$ of 
$C^{\infty}({\R }^2)$ into $C^0({\R }^2)$, where we have $v(x)=${\tiny $\left\{ \begin{array}{cl}\hspace{-7pt} 0, & \hspace{-8pt}
\mbox{if $x \le 0$} \\
\hspace{-5pt} h'(x), & \hspace{-8pt} \mbox{if $x >0$.} \end{array} \right. $} If $F \in 
C^{\infty}({\R }^2)$ then   
\begin{displaymath}
G(x,y) = X(F)(x,y) = \left\{ \begin{array}{cl} \frac{\partial F}{\partial x}(x,y), &\hspace{-8pt} \mbox{if $x \le 0$} \\
\rule{0pt}{12pt}\frac{\partial F}{\partial x}(x,y) + h'(x) \frac{\partial F}{\partial y}(x,y), & \hspace{-8pt} 
\mbox{if $x > 0$.} 
\end{array} \right. 
\end{displaymath}
So
\begin{displaymath}
\frac{{\partial }^{n+m}G}{\partial x^n \partial y^m}(x,y) = 
\left\{ \begin{array}{cl}
\hspace{-5pt}\frac{{\partial }^{n+m+1}F}{\partial x^{n+1}\partial y^m}(x,y), &\hspace{-7pt} \mbox{if $x \le 0$} \\
\hspace{-5pt} \frac{{\partial }^{n+m+1}F}{\partial x^{n+1}\partial y^m}(x,y) + 
\sum^n_{k=0} \binom{n}{k} h^{(n-k+1)}(x)\frac{{\partial }^{m+k+1}F}{\partial x^k \partial y^{m+1}}(x,y), 
&\hspace{-7pt} \mbox{if $x >0$.} 
\end{array} \right. 
\end{displaymath} 
For all $n \in {\Z }_{\ge 0}$ we have 
\begin{displaymath}
\lim_{x\rightarrow 0} \frac{{\partial}^{n+m}G}{\partial x^n\partial y^m} (x,0) 
= \frac{{\partial}^{n+m+1}F}{\partial x^n\partial y^m} (0, 0). 
\end{displaymath}
Since all the partial derivatives of $G$ are continuous at $(0,0)$, it follows that 
$G \in C^{\infty}({\R }^2)$. So $X$ is a derivation on $S$. Note that by definition 
$X(0,0) \ne (0,0)$. The derivation $X$ on $S$ is a vector field on $S_{\mathrm{reg}} = 
S \setminus \{ (0,0) \} $, which 
is a smooth manifold. \medskip 

\noindent We now show that a maximal integral curve of $X$ starting at $(0,0)$ has 
domain $\{ (0, (0,0)) \in \R \times S \}$.  
Consider the curve 
\begin{displaymath}
c : {\R }_{\ne 0} \rightarrow S_{\mathrm{reg}}: x \mapsto (x, y(x)) = (x,{\mathrm{e}}^{-(x^{-2})}).
\end{displaymath}  
Because $X$ is a vector field on $S_{\mathrm{reg}}$, we see that $c$ is 
the image of an integral curve of $X$ starting at $(1,{\mathrm{e}}^{-1})$. Now 
\begin{displaymath}
\gamma : {\R }_{\ne 0} \rightarrow S: t \mapsto 
(t+1, {\mathrm{e}}^{-(t+1)^{-2}}) = (x(t), y(t)) 
\end{displaymath}
satisfies $\dot{x}(t) = 1$ and $\dot{y}(t) = 2y(t){x(t)}^{-3}$, since 
\begin{align}
\frac{\dee y}{\dee t} (t) & = \frac{\dee }{\dee t}  {\mathrm{e}}^{-(t+1)^{-2}}  =  
2(t+1)^{-3}{\mathrm{e}}^{-(t+1)^{-2}} = 
2{x(t)}^{-3}y(t). \notag 
\end{align}%
Consequently, $\gamma $ is an integral curve of the vector field $Z = x^3 \frac{\partial }{\partial x} + 2y \frac{\partial }{\partial y}$ on ${\R }^2 \setminus \{ (0,0) \}$. The integral curves of $Z$ satisfy 
$\frac{\dee x}{\dee t} = x^3$ and $\frac{\dee y}{\dee t} = 2y$. Separating variables gives 
$\frac{\dee x}{x^3} = \dee t$ and $ \frac{\dee y}{2y} = \dee t$, 
which integrates to give the curve 
\begin{displaymath}
\gamma : (-\infty, 0] \rightarrow S_{\mathrm{reg}}: t \mapsto (x(t), y(t)) =  
\Big( \big(\frac{-1}{2t-1}\big)^{1/2}, {\mathrm{e}}^{(2t-1)} \Big) ,
\end{displaymath} 
that starts at $(1, {\mathrm{e}}^{-1})$ and has the same image as the 
curve $x \mapsto (x,{\mathrm{e}}^{-(x^{-2})})$ when $0< x \le 1$. Note that 
$\lim_{t \searrow -\infty} \gamma (t) = (0,0)$. So $\gamma $ is the maximal nonpositive time integral curve of the vector field $X$ on $S_{\mathrm{reg}}$ starting at $(1, {\mathrm{e}}^{-1})$. 
To see this suppose that the domain of a maximal nonpositive time 
integral curve $\Gamma $ of the derivation $X$ on $S$ starting at $(0,0)$ is an interval $(-\varepsilon, 0]$ for some $\varepsilon >0$. 
Then there is $t_0 \in (-\varepsilon, 0) \cap (-\infty, 0)$ such that $\Gamma (t_0) = \gamma (t_0) 
\in S_{\mathrm{reg}}$. By uniqueness of integral curves of vector fields $\Gamma = \gamma $ 
on $(-\varepsilon ,0)$. Hence $\Gamma $ is an extension of the maximal nonpositive 
time integral curve $\gamma $, which is a contradiction. 
So $\Gamma $ is only defined on $[0, \infty)$. 
Repeating the above argument for the integral curve ${\gamma }^{\vvee}: [0, \infty) \rightarrow S_{\mathrm{reg}}: t \mapsto (-1+t, {\mathrm{e}}^{-(-1+t)^{-2}})$ of the vector field $X$ starting at 
$(-1, {\mathrm{e}}^{-1})$ shows that ${\gamma }^{\vvee}$ is a maximal nonnegative time integral curve of $X$ because $\lim_{t \nearrow \infty}{\gamma }^{\vvee}(t) = (0,0)$. Hence $\Gamma $ is only defined at $t=0$. In other words, the nonzero derivation $X(0,0) = 
\frac{\partial }{\partial x}\rule[-5pt]{.5pt}{15pt}\raisebox{-4pt}{$\, \scriptscriptstyle (0,0)$}$ at 
the point $(0,0)$ is the lifted integral curve of the derivation $X$ on $S$ at $(0,0)$, which has a 
domain consisting only of $\{ (0, (0,0)) \in \R \times S \} $.  \hfill $\square $ \medskip   

\noindent \textbf{Theorem 3.2.3} Let $(S, C^{\infty}(S))$ be a subcartesian differential 
space and let $X$ be a derivation of $C^{\infty}(S)$. For every $x \in S$ there is 
a unique maximal lifted integral curve of $X$ starting at $x$. \medskip 

\noindent \textbf{Proof.} See \cite[thm 8]{cushman-sniatycki19}. \hfill $\square $ 

\subsection{Vector fields}

A vector field on a smooth manifold $M$ is not only a derivation of $C^{\infty}(M)$ but 
it has a unique integral curve starting at $m \in M$ and also it 
generates a local one-parameter group of local diffeomorphisms of $M$. On a 
subcartesian differential space $(S, C^{\infty}(S))$ not all derivations of $C^{\infty}(S)$ 
generate a local one-parameter group of diffeomorphisms of $S$, see example 3.3.2. 
We define a \emph{vector field on a subcartesian space} $(S, C^{\infty}(S))$ to 
be a derivation $X$ of $C^{\infty}(S)$ such that for every $x_0 \in S$ the following 
conditions hold: 1) there is a unique 
lifted integral curve of $X$ starting at $x_0$, and 2) there is an 
open neighborhood $U_{x_0}$ of $x_0$ and an ${\varepsilon }_{x_0} > 0$ such that 
for every $x \in U_{x_0}$ the interval $(-{\varepsilon }_{x_0}, {\varepsilon }_{x_0})$ 
is contained in the domain $I_x$ of the lifted integral curve 
${\gamma }_x: I_x \subseteq \R \rightarrow TS$ starting at $x$ of the derivation $X$ and 
the map 
\begin{displaymath}
{\mathrm{e}}^{tX}: U_{x_0} \subseteq S \rightarrow S: x \mapsto ({\tau }_S \comp {\gamma }_x)(t)
\end{displaymath}
is defined for every $t \in (-{\varepsilon }_{x_0}, {\varepsilon }_{x_0})$ and is a 
diffeomorphism of $U_{x_0}$ onto the open subset ${\mathrm{e}}^{tX}(U_{x_0})$ of 
$S$. \medskip 

Recall that a subset $A$ of a topological space $\mathcal{T}$ is \emph{locally closed}
if for every point $x \in A$ there is an open set $U_x \in \mathcal{T}$ containing $x$ such that 
$A\cap U_x$ is closed in $U_x$. For a locally closed differential space we have a much simpler criterion for recognizing whether a derivation is a vector field. \medskip 

\noindent \textbf{Proposition 3.3.1} Let $S$ be a locally closed subcartesian differential 
space. A derivation $X$ of $C^{\infty}(S)$ with unique integral curves 
is a vector field on $S$ if the domain of 
every maximal integral curve of $X$ is an open interval in $\R $ which contains $0$. \medskip 

\noindent \textbf{Proof.} See \cite[prop 3.2.3, p.35]{sniatycki}. \hfill $\square $ \medskip 

The hypothesis that $S$ is locally closed is necessary as the following example shows. \medskip 

\par \noindent \textbf{Example 3.3.2} Consider the subset of ${\R }^2$ given by 
\begin{displaymath}
S = \{ x = (x_1,x_2) \in {\R }^2 \setrule \, x^2_1 +(1-x_2)^2 < 1 \, \, \mathrm{or} \, \, 
x_2 = 0 \} . 
\end{displaymath}
Then $(S, C^{\infty}({\R }^2)|_S)$ is a subcartesian differential space 
but $S$ is not locally closed at $(0,0)$. The vector field $Y = \frac{\partial }{\partial x_1}$ 
on ${\R }^2$ restricts to a derivation $X$ of $C^{\infty}(S)$. For every $x \in {\R }^2$ 
we have ${\mathrm{e}}^{tY}x = (x_1+t, x_2)$ for all $t \in \R $. All the maximal integral 
curves of $X$ have open domains which contain $0$. But ${\mathrm{e}}^{tX}$ is 
\emph{not} a local one-parameter group of local diffeomorphisms of $S$. 
\hfill $\square $ 

\section{Differentiable structure}

For each $n \in {\Z }_{\ge 0}$ let ${\mathfrak{R}}_n$ be the set of all subsets of which are 
open in ${\R }^n$. Let $\mathfrak{R} = \bigcup_{n \in {\Z }_{\ge 0}} {\mathfrak{R}}_n$. Let $S$ be a topological space with 
$\mathcal{F}$ a subset of $C^0(S)$, the set of continuous functions on $S$. 
For every $U \in \mathfrak{R}$ let $C^0(U,S)$ be the set of continuous mappings 
of $U$ into $S$. Let
\begin{equation}
\mathcal{F}(U) = \{ g \in C^0(U,S) \setrule \, g^{\ast }f = f \comp g \in C^{\infty}(U)\, \, 
\mbox{for every $f \in \mathcal{F}$} \}. 
\label{eq-s3nwone}
\end{equation}
We say that $\mathcal{F}$ is a \emph{differentiable structure} on $S$ if, given 
$f \in C^0(S)$ such that $g^{\ast }f \in C^{\infty}(U)$ for every $U \in \mathfrak{R}$ 
and every $g \in \mathcal{F}(U)$, then $f \in \mathcal{F}$. A \emph{differentiable space} 
is a pair $\mathbf{S} = (S, \mathcal{F})$, where $\mathcal{F} ={\mathcal{F}}_{\mathbf{S}}$ is a differentiable structure on $S$. \medskip 

\noindent \textbf{Proposition 4.1} Given a topological space $S$ and a subset $\mathcal{F}$ of 
$C^0(S)$, there is a unique differentiable structure ${\mathcal{F}}^{\ast }$ on $S$ such that 
${\mathcal{F}}^{\ast }(U) = \mathcal{F}(U)$ for every $U \in \mathfrak{R}$. \medskip 

\noindent \textbf{Proof.} Let 
\begin{displaymath}
{\mathcal{F}}^{\ast } = \{ f \in C^0(U) \setrule \, f \comp g \in C^{\infty}(U ) \, \, 
\mbox{for every $U \in \mathfrak{R}$ and every $g \in \mathcal{F}(U)$} \} . 
\end{displaymath}
If $f \in \mathcal{F}$, then for every $U \in \mathfrak{R}$ and 
every $g \in \mathcal{F}(U)$ we have $g^{\ast }f \in C^{\infty}(U)$. So $f \in 
{\mathcal{F}}^{\ast }$. Thus $\mathcal{F} \subseteq {\mathcal{F}}^{\ast }$. Suppose 
that $U \in \mathfrak{R}$ and $g \in \mathcal{F}(U)$. Then $g^{\ast }f \in C^{\infty}(U)$ 
for every $f \in \mathcal{F}$. Hence $g \in {\mathcal{F}}^{\ast }(U)$. 
So $\mathcal{F}(U) \subseteq {\mathcal{F}}^{\ast }(U)$. Let $g \in {\mathcal{F}}^{\ast }(U)$. 
Then $f \comp g \in C^{\infty}(U)$ for every $f \in {\mathcal{F}}^{\ast }$. Since 
$\mathcal{F} \subseteq {\mathcal{F}}^{\ast }$, it follows that $f \comp g \in C^{\infty}(U)$ for 
every $f \in \mathcal{F}$. Hence $g \in \mathcal{F}(U)$. So ${\mathcal{F}}^{\ast }(U) 
\subseteq \mathcal{F}(U)$. Consequently, ${\mathcal{F}}^{\ast }(U) = 
\mathcal{F}(U) $ for every $U \in \mathfrak{R}$. \hfill $\square $ \medskip 

\noindent ${\mathcal{F}}^{\ast }$ is the differentiable structure on $S$ \emph{generated} by 
$\mathcal{F}$. Clearly ${\mathcal{F}}^{\ast }$ is a subring of $C^0(S)$ with unit, the constant function $\mathbf{1}:S \rightarrow \R :x \mapsto 1$.  \medskip 

\noindent \textbf{Proposition 4.2} Let $M$ be a smooth manifold. Then $C^{\infty}(M)$ is a differentiable structure on $M$.  \medskip 

\noindent \textbf{Proof.} Let $m \in M$ and let $V$ be an open subset of $M$ containing 
$m$ such that $\varphi : V \subseteq M \rightarrow \varphi (V) \subseteq {\R }^n$ is a 
diffeomorphism. Suppose that $U \in \mathfrak{R}$ and $g: U \rightarrow V \subseteq M$ is 
a continuous mapping such that $m = g(u)$ for some $u \in U$. For some open subset 
$W$ of $m$ in $V$ we have $U' = g^{-1}(W)$ is an open subset of $U$ containing the 
point $u$. \medskip 

Let $\mathcal{F} = C^{\infty}(M)$. We now show that $\mathcal{F}(U') = C^{\infty}(U', U)$. 
Suppose that $g: U' \rightarrow U \subseteq M$ satisfies the condition 
$f \comp g \in C^{\infty}(U')$ for every $f \in C^{\infty}(M)$, that is, $g \in \mathcal{F}(U')$. 
If $g \in C^{\infty}(U', U)$, then $f \comp g \in C^{\infty}(U')$ for every $f \in C^{\infty}(M)$. 
So $g \in \mathcal{F}(U')$. Thus $C^{\infty}(U', U) \subseteq \mathcal{F}(U')$. We now show 
that $\mathcal{F}(U') \subseteq C^{\infty}(U', U)$. Suppose that $g \in \mathcal{F}(U')$. Then 
\begin{displaymath}
f \comp g = (f \comp {\varphi }^{-1}) \comp (\varphi \comp g) = 
(\varphi \comp g)^{\ast } \big( f \comp {\varphi }^{-1}). 
\end{displaymath}
Since $f \in C^{\infty}(M) = \mathcal{F}$, it follows that the function 
$f \comp {\varphi }^{-1}: \varphi (V) \subseteq {\R }^n \rightarrow \R$ lies in 
$C^{\infty}(\varphi (V))$. So $(\varphi \comp g)^{\ast } \big( C^{\infty}(\varphi (V)) \big) 
\subseteq C^{\infty}(U')$. Thus the mapping $\varphi \comp g: U'  \rightarrow \varphi (V) 
\subseteq {\R }^n$ is smooth. Hence the mapping $g = {\varphi }^{-1}\comp (\varphi \comp g): 
U' \rightarrow U \subseteq M$ is smooth. This shows that $\mathcal{F}(U') \subseteq 
C^{\infty}(U', U)$ as desired. So $\mathcal{F}(U') =  C^{\infty}(U', U)$. \medskip

Suppose that $f \in C^0(M)$ and $f \comp g \in C^{\infty}(U')$ for every $U' \in \mathfrak{R}$ 
and every $g \in \mathcal{F}(U')$. Suppose that $f \notin C^{\infty}(M)$. Then for some 
$m \in M$ there is an open neighborhood $U$ of $m$ in $M$ such that 
$(U, \varphi )$ is a chart with $\varphi : U \subseteq M \rightarrow \varphi (U) \subseteq {\R }^n$ 
a diffeomorphism and the function $f \comp {\varphi }^{-1}: \varphi (U) \subseteq {\R }^n 
\rightarrow \R $ is not smooth at $\varphi (m)$. Set $U' = \varphi (U)$ and let 
$g = {\varphi }^{-1}: U' = \varphi (U) \subseteq {\R }^n \rightarrow U \subseteq M$. Then 
$g \in C^{\infty}(U', U) =  \mathcal{F}(U')$ but $g^{\ast }f: U' \rightarrow \R $ is not 
smooth at $\varphi (m) \in U'$. This contradicts our hypothesis that $g \in \mathcal{F}(U')$. 
Hence $f \in C^{\infty}(M)$. In other words, $C^{\infty}(M)$ is a differentiable structure on $M$. 
\hfill $\square $ \medskip

Let $S$ be a differential space. Define
\begin{equation}
{\Gamma }_0C^{\infty}(S) = \{ c \in C^0(\R ,S) \setrule \, f \comp c \in C^{\infty}(\R ) \, \, 
\mbox{for all $f \in C^{\infty}(S)$} \} .
\label{eq-s3ss3threevnw}
\end{equation}
From the definition of the differential space topology on $S$, it follows that 
${\Gamma }_0C^{\infty}(S) = C^0(\R ,S)$. Also define
\begin{equation}
{\Phi }_0{\Gamma }_0C^{\infty}(S) = \{ f \in C^0(S) \setrule \, 
f \comp c \in C^{\infty}(\R ) \, \, \mbox{for all $c \in {\Gamma }_0C^{\infty}(S)$} \} . 
\label{eq-s3ss3fourvnw}
\end{equation}
Clearly $C^{\infty}(S) \subseteq {\Phi }_0{\Gamma }_0C^{\infty}(S)$. We say that 
$C^{\infty}(S)$ is \emph{continuously reflexive} if ${\Phi }_0{\Gamma }_0C^{\infty}(S) \subseteq 
C^{\infty}(S)$. We have \medskip

\noindent \textbf{Theorem 4.3} If $C^{\infty}(S)$ is a continuously reflexive differential 
structure on the differential space $(S, C^{\infty}(S))$, then $C^{\infty}(S)$ is a differentable 
structure on $S$ in the sense of Smith. \medskip 

\noindent \textbf{Proof.} Let $S$ be a topological space 
with $\mathcal{F}$ a subset of $C^0(S)$, the set of continuous functions on $S$. 
For every $U \in \mathfrak{R}$ let $C^0(U,S)$ be the set of continuous mappings 
of $U$ into $S$. Set 
\begin{equation}
\mathcal{F}(U) = \{ g \in C^0(U,S) \setrule \, g^{\ast }f = f \comp g \in C^{\infty}(U) \, \, 
\mbox{for every $f \in \mathcal{F}$} \}. 
\label{eq-3nwtwodagger}
\end{equation}
Recall that $\mathcal{F}$ is a differentiable structure on $S$ if $f \in C^0(S)$ lies 
in $\mathcal{F}$ whenever the condition 
\begin{equation}
g^{\ast }f \in C^{\infty}(U) \, \, \mbox{for every $U \in \mathfrak{R}$ and every $g \in \mathcal{F}(U)$}  
\label{eq-3nwfour}
\end{equation}
is satisfied. Take $\mathcal{F} = C^{\infty}(S)$. Then definition (\ref{eq-3nwtwodagger}) reads: for every $U \in \mathfrak{R}$ 
\begin{equation}
C^{\infty}(S)(U) = \{ g \in C^0(U,S) \setrule \, g^{\ast }f \in C^{\infty}(U) \, \, 
\mbox{for every $f \in C^{\infty}(S)$} \} . 
\label{eq-3nwfourstar}
\end{equation}
So condition (\ref{eq-3nwfour}) that $C^{\infty}(S)$ is a differentiable structure on $S$ reads: a continuous function $f$ on $S$ lies in $C^{\infty}(S)$ whenever the condition 
\begin{equation}
g^{\ast }f \in C^{\infty}(U)\, \, \mbox{for every $U \in \mathfrak{R}$ and every $g \in C^{\infty}(S)(U)$}
\label{eq-3nwfourdoublestar}
\end{equation} 
is satisfied. Let $f \in C^0(S)$ and suppose that condition (\ref{eq-3nwfourdoublestar}) and definition 
(\ref{eq-3nwfourstar}) hold. Specializing to the case when $U = \R \in {\mathfrak{R}}_1$ definition 
(\ref{eq-3nwfourstar}) reads: 
\begin{align}
C^{\infty}(S)(\R ) & = \{ g \in C^0(\R ,S) \setrule \, g^{\ast }f \in C^{\infty}(\R) \, \, 
\mbox{for every $f \in C^{\infty}(S)\}  = {\Gamma }_0C^{\infty}(S)$}  ,  \notag 
\end{align}
using (\ref{eq-s3ss3threevnw}), and condition (\ref{eq-3nwfourdoublestar}) reads: 
\begin{equation}
g^{\ast}f \in C^{\infty}(\R ) \, \, \mbox{for every $g \in {\Gamma }_0C^{\infty}(S)$.} 
\label{eq-3nwfivenew}
\end{equation}
From equation (\ref{eq-s3ss3fourvnw}) it follows that $f \in {\Phi }_0\big( {\Gamma }_0(C^{\infty}(S)) \big)$. 
By hypothesis $C^{\infty}(S)$ is continuously reflexive, that is, 
${\Phi }_0\big( {\Gamma }_0(C^{\infty}(S)) \big) = C^{\infty}(S)$. So $f \in C^{\infty}(S)$. 
Thus $C^{\infty}(S)$ is a differentiable structure on $S$. \hfill $\square $ \medskip 

Let $\mathbf{S} = (S, {\mathcal{F} }_{ \mathbf{S} } )$ and ${\mathbf{S}}' = 
(S', {\mathcal{F}}_{ {\mathbf{S}}' } )$ be differentiable 
spaces and let $h:\mathbf{S} \rightarrow {\mathbf{S} }'$ be a continuous map. $h$ is a 
\emph{differentiable mapping} if and only if $h^{\ast }({\mathcal{F}}_{\mathbf{S}' }) \subseteq 
{\mathcal{F}}_{\mathbf{S}}$.  \medskip 

\noindent \textbf{Proposition 4.4} A continuous function $f:S \rightarrow \R $ is a differentiable map $\mathbf{f}:\mathbf{S} \rightarrow \mathbf{R} = (\R , C^{\infty}(\R ))$ if and only if 
$f \in {\mathcal{F}}_{\mathbf{S}}$. \label{pg-ten} \medskip 

\noindent \textbf{Proof:} This follows immediately from the definition of smooth map 
between differentiable spaces. \hfill $\square $ \medskip

\noindent \textbf{Proposition 4.5} A map $\varphi : (S_1, C^{\infty}(S_1)) \rightarrow 
(S_1, C^{\infty}(S_1))$ of differential spaces with continuously reflexive differential structures 
is smooth if and only if it is a differentable 
mapping between the differentiable spaces ${\mathbf{S}}_1 = (S_1, C^{\infty}(S_1))$ and 
${\mathbf{S}}_2 = (S_2, C^{\infty}(S_2))$. \medskip 

\noindent \textbf{Proof.} The map $\varphi $ between the differential spaces 
$S_1$ and $S_2$ is smooth if and only if ${\varphi }^{\ast } \big( C^{\infty}(S_2) \big) 
\subseteq C^{\infty}(S_1)$ if and only if it is a differentable map between differentiable 
spaces ${\mathbf{S}}_1$ and ${\mathbf{S}}_2$. \hfill $\square $ \medskip  

We now give some results on differentable structures that we will need when discussing 
differential forms on differentiable spaces. \medskip 

\noindent \textbf{Proposition 4.6} If $g_1, \ldots , g_n \in {\mathcal{F}}_{\mathbf{S}}$,  
$\phi \in C^{\infty}({\R }^n)$, and $U \in \mathfrak{R}$, then 
\par \noindent \hspace{-.2in}\begin{tabular}{l}
\vspace{-.15in} \\
\hspace{.45in} $f = \phi (g_1, \ldots , g_n): U \rightarrow \R : x \mapsto f(x) = 
\phi \big( g_1(x), \ldots , g_n(x) \big) $ \\
\end{tabular}
\par \noindent lies in ${\mathcal{F}}_{\mathbf{S}}$. \medskip 

\noindent \textbf{Proof.} For every $U \in \mathfrak{R}$ and every $h: U \rightarrow S$ in 
${\mathcal{F}}_{\mathbf{S}}(U)$ we have $g_i \comp h \in C^{\infty}(U)$. Hence 
$f \comp h \in C^{\infty}(U)$, which implies that $f \in {\mathcal{F}}_{\mathbf{S}}$. 
\hfill $\square $ \medskip 

Let ${\mathbf{S}}' = (S', {\mathcal{F}}_{{\mathbf{S}}'})$ and $({\mathbf{S}}'', 
{\mathcal{F}}_{{\mathbf{S}}''})$ be 
differentiable spaces. Define the product ${\mathbf{S}}' \times {\mathbf{S}}'' = 
(S' \times S'', {\mathcal{F}}_{{\mathbf{S}}' \times {\mathbf{S}}''})$, where 
${\mathcal{F}}_{{\mathbf{S}}' \times {\mathbf{S}}''}$ is the set of $f \in C^0(S' \times S'' )$ such that 
$f (g', g'') \in C^{\infty}(U)$ for every $U \in \mathfrak{R}$ and every 
$g' \in {\mathcal{F}}_{{\mathbf{S}}'}$ and $g'' \in {\mathcal{F}}_{{\mathbf{S}}''}$. \medskip 

\noindent \textbf{Lemma 4.7} ${\mathcal{F}}_{{\mathbf{S}}' \times {\mathbf{S}}''}$ is a differentiable 
structure.  \medskip 

\noindent \textbf{Proof.} Let $f: S \rightarrow \R $ be a continuous function such that $f \comp h 
\in C^{\infty}(U)$ for every $U \in \mathfrak{R}$ and every 
$h \in {\mathcal{F}}_{{\mathbf{S}}' \times {\mathbf{S}}''}$. Consider 
$h = (h', h'')$ with $h' \in {\mathcal{F}}_{{\mathbf{S}}'}(U)$ and 
$h'' \in {\mathcal{F}}_{{\mathbf{S}}''}(U)$. 
For every $g \in {\mathcal{F}}_{{\mathbf{S}}' \times {\mathbf{S}}''}$ we have 
$g \comp h \in C^{\infty}(U)$. 
This implies that $h \in {\mathcal{F}}_{{\mathbf{S}}' \times {\mathbf{S}}''}(U)$. Thus 
$f \comp h \in C^{\infty}(U)$, that is, $f \in {\mathcal{F}}_{{\mathbf{S}}' \times {\mathbf{S}}''}$. 
\hfill $\square $ \medskip 

\noindent  Hence the product ${\mathbf{S}}' \times {\mathbf{S}}'' = (S' \times S'', 
{\mathcal{F}}_{{\mathbf{S}}' \times {\mathbf{S}}''})$ is a differentiable space. \medskip 

Suppose that $g'_1, \ldots , g'_n \in {\mathcal{F}}_{{\mathbf{S}}'}$ and 
$g''_1, \ldots , g''_m \in {\mathcal{F}}_{{\mathbf{S}}''}$ are given and that 
$\phi \in C^{\infty}( {\R }^{n+m} )$. 
Let ${\pi }': S' \times S'' \rightarrow S'$ and ${\pi }'': S' \times S'' \rightarrow S''$ be the natural 
projections. Let $f'_i = g'_i \comp {\pi }'$ and $f''_j = g''_j \comp {\pi }''$. Then 
\begin{equation}
f = \phi (f'_1, \ldots , f'_n, f''_1, \ldots , f''_m) 
\label{eq-s3nwzero}
\end{equation}
lies in ${\mathcal{F}}_{{\mathbf{S}}' \times {\mathbf{S}}''}$. Let 
${\overline{\mathcal{F}}}_{{\mathbf{S}}' \times {\mathbf{S}}''}$ be the subset of all functions in 
${\mathcal{F}}_{{\mathbf{S}}' \times {\mathbf{S}}''}$ which can be written in the form 
(\ref{eq-s3nwzero}). Then ${\overline{\mathcal{F}}}_{{\mathbf{S}}' \times {\mathbf{S}}''}$ is a linear subspace of the real vector space ${\mathcal{F}}_{{\mathbf{S}}' \times {\mathbf{S}}''}$. \medskip 

 A differentiable map $h: {\mathbf{S}}' \times {\mathbf{S}}'' \rightarrow \mathbf{T}$ is 
\emph{well behaved} if $g \comp h \in {\overline{\mathcal{F}}}_{{\mathbf{S}}' \times {\mathbf{S}}''}$ for 
every $g \in {\mathcal{F}}_{\mathbf{T}}$. Clearly the natual projections ${\pi }'$ and ${\pi }''$ are well behaved differentiable mappings. We will need this notion in section 5, when proving the Poincar\'{e} lemma. \medskip 

The collection of differentiable spaces and differentiable maps form a category $\mathcal{D}$ under composition. The category of smooth manifolds and smooth maps 
is a full subcategory of $\mathcal{D}$. Let $\mathcal{K}$ be the category 
of rings with unit and ring homomorphisms. We have a contravariant functor 
$F: \mathcal{D} \rightarrow \mathcal{K}$, which 
sends a differentiable space $\mathbf{S}$ to its differentiable structure 
${\mathcal{F}}_{\mathbf{S}}$, 
which is a ring with unit and sends a differentiable map between differentiable spaces to 
a ring homomorphism between the differential structures ${\mathcal{F}}_{{\mathbf{S}}'}$ and 
${\mathcal{F}}_{\mathbf{S}}$. More precisely, if $h: \mathbf{S} = (S, {\mathcal{F}}_{\mathbf{S}}) \rightarrow {\mathbf{S}}' = (S', {\mathcal{F}}_{{\mathbf{S}}'})$ is a differentiable map, then $h^{\ast }: {\mathcal{F}}_{{\mathbf{S}}'} \rightarrow {\mathcal{F}}_{\mathbf{S}}: f \mapsto f \comp h$ is a ring homomorphism. Let $U$ be an open subset of $S$ and let ${\mathcal{F}}_U$ be the differentiable structure on $U$ generated by $f|_U$ where $f \in {\mathcal{F}}_{\mathbf{S}}$. Let ${\mathbf{S}}|_U = (U, {\mathcal{F}}_U)$ and let ${\iota }_{\mathbf{S}}|_U: \mathbf{S}|_U \rightarrow \mathbf{S}$ be the inclusion map. This defines a \emph{local} category $\mathcal{D}$. Similarly, we can define 
a local category ${\mathcal{S}}^0$ with objects consisting of pairs of an open subset $U$ of the topological space $S$ and the set of continuous functions on $S$ restricted to $U$. 
Thus the objects of ${\mathcal{S}}^0$ are $S_U = (U, C^0(U))$. The maps in ${\mathcal{S}}^0$ are 
the inclusions ${\iota }_U: U \mapsto S$. For every object 
$\mathbf{S}$ in the local category $\mathcal{D}$ we have a covariant functor $T_{\mathbf{S}}: {\mathcal{S}}^0 \rightarrow \mathcal{D}$ whose value at ${\mathbf{S}}_U$ is the open subset $U$ of $S$ and whose value at the map ${\iota }_{U}|_V$ is the map 
$({\iota }_{\mathbf{S}|_U})|_V$ for every open subset $U$ and $V$ of $S$ such that $V \subseteq U$. For every $\mathbf{S} \in \mathcal{D}$ the map 
\begin{displaymath}
F^0 = F \comp T_{\mathbf{S}}: {\mathcal{S}}^0 \rightarrow \mathcal{K}: 
U \subseteq S \mapsto {\mathcal{F}}_{{\mathbf{S}}|_U}
\end{displaymath}
is a \emph{sheaf} on $S$ with values in the category $\mathcal{K}$. 

\section{Exterior algebra on a differentiable space}

In this section we construct an exterior algebra of differential forms on a differentiable space 
and prove a version of de Rham's theorem. We follow the arguments of Smith \cite{smith} quite closely and give the details. 

\subsection{Differential forms}

In this subsection we construct the graded exterior algebra of differential forms and 
exterior derivative on the differentiable space $S$. \medskip 

A \emph{$p$-dimensional cube} $J^p \subseteq {\R }^p$ with $p >0$ is a closed subset of 
${\R }^p$ bounded by $2p$ axis parallel hyperplanes. A map $\sigma : J^p \subseteq {\R }^p \rightarrow S$ is a \emph{singular $p$-cube} in the differentiable space $\mathbf{S} = 
(S, {\mathcal{F}}_{\mathbf{S}})$ 
if there is a $U\in {\mathfrak{R}}_p$ with $J^p \subseteq U \subseteq {\R }^p$ and a map 
$f: U \rightarrow S$ in ${\mathcal{F}}_{\mathbf{S}}(U)$ such that 
$\sigma = f|_{J^p}$. We say that the singular $p$-cube $\sigma $ is \emph{represented by} the 
mapping $f \in {\mathcal{F}}_{\mathbf{S}}(U)$. Let $K_p(\mathbf{S})$ be the set of all singular $p$-cubes in $\mathbf{S}$. Let $C^p(\mathbf{S})$ be the set of all real valued functions $\alpha : K_p(\mathbf{S}) \rightarrow \R $. 
An element of $C^p(\mathbf{S})$ is a \emph{singular $p$-cochain} on $\mathbf{S}$. Clearly, 
$C^p(\mathbf{S})$ is a real vector space. When $p=0$ let $J^0 = \{ 0 \} $ in ${\R }^0 = \{ 0 \} $. A singular $0$-cube is a mapping $\sigma : J^0 \subseteq {\R }^0 \rightarrow S$ which sends 
$0$ to a point $x \in S$. Hence $K_0(\mathbf{S})$, the set of singular $0$-cubes in 
$\mathbf{S}$, is $S$. The set $C^0(\mathbf{S})$ of singular $0$-cochains on $\mathbf{S}$ is 
${\mathcal{F}}_{\mathbf{S}}$, the differentiable structure of $\mathbf{S}$. \medskip 

For every $p > 0$ define the mapping
\begin{equation}
{\lambda }_p: {\mathcal{F}}^{p+1}_{\mathbf{S}} = \overbrace{{\mathcal{F}}_{\mathbf{S}} 
\times \cdots \times {\mathcal{F}}_{\mathbf{S}}}^{p+1} \rightarrow C^p(\mathbf{S}): 
(f_0, \ldots , f_p) \mapsto {\lambda }_p(f_0, \ldots , f_p),  
\label{eq-s4onenw}
\end{equation}
where $f_i \in {\mathcal{F}}_{\mathbf{S}}$. To specify ${\lambda }_p(f_0, \ldots , f_p)$, we need to 
give its value on a singular $p$-cube $\sigma : J^p \rightarrow S$ represented by $f \in 
{\mathcal{F}}_{\mathbf{S}}(U)$. Define
\begin{equation}
{\lambda }_p(f_0, \ldots , f_p)(\sigma ) = 
\int_{J^p} g_0 \det DG (t_1, \ldots ,t_p) \, \dee t_1 \cdots \dee t_p, 
\label{eq-s4twonw}
\end{equation}
where 
\begin{displaymath}
\begin{array}{l}
G: U \subseteq {\R }^p \rightarrow {\R }^p : \\
\hspace{.25in} \mathbf{t} = (t_1, \ldots , t_p) \mapsto \big( g_1(\mathbf{t}), \ldots , g_p(\mathbf{t}) \big) 
= \big( (f_1 \comp f)(\mathbf{t}), \ldots , (f_p \comp f)(\mathbf{t}) \big) 
\end{array}
\end{displaymath}
is a smooth mapping and $g_0 = f_0 \comp f$. The Jacobian $\det DG$ of the mapping $G$ will be denoted by $\frac{\partial ( g_1, \ldots , g_p)}{\partial (t_1, \ldots , t_p)}$. 
Since the right hand side of (\ref{eq-s4twonw}) depends only on the singular 
$p$-cube $\sigma $ and not on its representation $f$, it follows that 
${\lambda }_p(f_0, \ldots , f_p)$ is a singular $p$-cochain. Hence the mapping 
${\lambda }_p$ (\ref{eq-s4onenw}) is well defined. When $p=0$ the set $C^0(\mathbf{S})$ of 
singular $0$-cochains on $\mathbf{S}$ is ${\mathcal{F}}_{\mathbf{S}}$. The map 
${\lambda }_0: {\mathcal{F}}_{\mathbf{S}} \rightarrow C^0(\mathbf{S}) = {\mathcal{F}}_{\mathbf{S}}$ is 
the identity map, because ${\lambda }_0(f_0)\sigma = f_0 \comp \sigma $, where 
$\sigma $ is a singular $0$-cube. \medskip 

For $p > 0$ let $G^p(\mathbf{S})$ be the real linear vector subspace of $C^p(\mathbf{S})$ spanned 
by the image of the map ${\lambda }_p$. When $p=0$ set $G^0(\mathbf{S}) = 
{\mathcal{F}}_{\mathbf{S}}$. We call $G^p(\mathbf{S})$ the real vector space of 
\emph{differential $p$ forms} on the differentiable space $\mathbf{S}$. \medskip 

Let $F:\mathbf{S} \rightarrow {\mathbf{S}}'$ be a differentiable mapping. Then $F$ induces a real 
linear mapping 
\begin{equation}
F^{\ast }:C^p({\mathbf{S}}') \rightarrow C^p(\mathbf{S}): 
{\lambda }_p(f'_1, \ldots , f'_p) \mapsto {\lambda }_p(f_1, \ldots , f_p), 
\label{eq-s4threenw}
\end{equation}
where $f_i = f'_i \comp F$. We now show that the induced map $F^{\ast }$ on cochains is 
well defined. Let $\sigma :J^p \rightarrow S$ be a singular $p$-cube represented by 
$f \in {\mathcal{F}}_{\mathbf{S}}(U)$. Then $F \comp \sigma \in 
C^p({\mathbf{S}}')$. Check: $F \comp \sigma : J^p \rightarrow S'$ and $F \comp f: U \rightarrow S'$ 
is in ${\mathcal{F}}_{{\mathbf{S}}'}(U)$ because $g \comp (F \comp f) \in C^{\infty}(U)$ 
for every $g \in {\mathcal{F}}_{{\mathbf{S}}'}$. The last assertion follows \linebreak 
because $g \comp (F \comp f) = 
F^{\ast }(g) \comp f \in C^{\infty}(U)$ for $F^{\ast }(g) \in {\mathcal{F}}_{\mathbf{S}}$, 
since $f \in {\mathcal{F}}_{\mathbf{S}}(U)$. Now $(F \comp f)|_{J^p} = F \comp (f|_{J^p}) = 
F \comp \sigma $. Thus $F \comp \sigma \in C^p({\mathbf{S}}')$. 
The following computation finishes the verification of equation (\ref{eq-s4threenw}). 
\begin{align}
{\lambda }_p(f'_1, \ldots , f'_p)\big( F \comp \sigma \big) & = 
\int_{J^p} (f'_0 \comp F) \frac{\partial \big( f'_1 \comp (F \comp f), \ldots , f'_p \comp (F \comp f) \big)}
{\partial (t_1, \ldots , t_p)} \, \dee t_1 \cdots \dee t_p \notag \\
& = \int_{J^p} f_0 \frac{\partial \big( f_1 \comp f, \ldots ,f_p \comp f \big)}{\partial (t_1, \ldots , t_p) } \, 
\dee t_1 \cdots \dee t_p \notag \\
& = {\lambda }_p(f_0, \ldots , f_p)(\sigma ). \tag*{$\square $}
\end{align} 

\noindent \textbf{Example 5.1.1} Let $U \in \mathfrak{R}$. Then $U $ is a smooth manifold. 
So $\mathbf{U} = (U , C^{\infty}(U))$ is a differentiable space with differential 
structure $C^{\infty}(U)$. Let $F^p(\mathbf{U})$ be the real vector space of 
differentiable $p$-forms on $U$. The map 
\begin{equation}
\Phi : G^p(\mathbf{U}) \rightarrow F^p(\mathbf{U}): 
{\lambda }_p(g_0, \ldots , g_p) \mapsto g_0 \dee g_1 \wedge \cdots \wedge \dee g_p , 
\label{eq-s4fournw}
\end{equation}
where $g_i \in C^{\infty}(U)$, is an isomorphism. 
The mapping $\Phi $ is well defined. For every singular $p$-cube 
$\sigma :J^p \subseteq  U \rightarrow U $ in $U$ we have 
\begin{equation}
\int_{J^p} g_0 \dee g_1 \wedge \cdots \wedge \dee g_p = {\lambda }_p(g_0, \ldots , g_p)\sigma . 
\label{eq-s4fivenw}
\end{equation}
Check: The singular $p$-cube $\sigma $ is represented by $\mathrm{id}: U \rightarrow U$. By definition 
\begin{align}
{\lambda}_p (g_0, \ldots , g_p)\sigma & = \int_{J^p} (g_0 \comp \mathrm{id}) \,  
\frac{\partial \big( g_1 \comp \mathrm{id} , \ldots , g_p \comp \mathrm{id} \big) }
{\partial (t_1, \ldots , t_p)} \dee t_1 \cdots \dee t_p \notag \\
& = \int_{J^p} g_0 \frac{\partial \big( g_1, \ldots , g_p \big) }{\partial (t_1, \ldots , t_p)} \dee t_1 \cdots \dee t_p 
\notag \\
& = \int_{J^p} g_0 \dee g_1 \wedge \cdots \wedge \dee g_p . \notag 
\end{align}
Clearly the map $\Phi $ (\ref{eq-s4fournw}) is a real linear mapping. It is injective, for if 
${\lambda }_p(g_0, \ldots , g_p) $ $=0$, then ${\lambda }_p(g_0, \ldots g_p)\sigma =0$ for 
all singular $p$-cubes $\sigma $. From (\ref{eq-s4fivenw}) we obtain $\int_{J^p} g_0 \dee g_1 \wedge 
\cdots \wedge \dee g_p = 0$ for every $p$-cube $J^p$ in ${\R }^p$. Thus $g_0 \dee g_1 \wedge 
\cdots \wedge \dee g_p =0$. Equation (\ref{eq-s4fivenw}) shows that the mapping $\Phi $ is 
surjective. Consequently, $\Phi $ is an isomorphism. \hfill $\square $ \medskip  

If $\alpha = {\lambda }_p(g_0, \ldots , g_p) \in G^p(\mathbf{S})$ with $g_i \in 
{\mathcal{F}}_{\mathbf{S}}$ 
and if $f \in {\mathcal{F}}_{\mathbf{S}}(U)$ represents the singular $p$-cube $\sigma : J^p  \rightarrow S$, then $f^{\ast }(\alpha ) = (g_0 \comp f) \dee \, (g_1\comp f) \wedge \cdots \wedge 
\dee \, (g_p \comp f) \in F^p(\mathbf{U})$, where $g_i \comp f \in C^{\infty}(U)$.

\subsection{Exterior derivative}

First we prove the existence of an \emph{exterior product} on the space of differential forms on the differentiable space $\mathbf{S}$. \medskip

\noindent \textbf{Proposition 5.2.1} There is a unique bilinear map 
\begin{displaymath}
\begin{array}{l}
\Lambda : G^p(\mathbf{S}) \times G^q(\mathbf{S}) \rightarrow G^{p+q}(\mathbf{S}): \\
\hspace{.25in} \big( {\lambda }_p(f_0, \ldots , f_p) , {\lambda }_q({\overline{f}}_0, \ldots , 
{\overline{f}}_q ) \big) \mapsto {\lambda }_{p+q}\big( f_0, {\overline{f}}_0, f_1, \ldots , f_p, 
{\overline{f}}_1, \ldots , {\overline{f}}_q \big) , 
\end{array}
\end{displaymath}
where $f_i$, ${\overline{f}}_j \in G^0(\mathbf{S}) = {\mathcal{F}}_{\mathbf{S}}$. \medskip 

\noindent \textbf{Proof.} Let $\alpha = {\lambda }_p(f_0, \ldots , f_p)$ and 
$\beta = {\lambda }_q({\overline{f}}_0, \ldots , {\overline{f}}_q)$. Let $\sigma :J^{p+q} \rightarrow S$ be a singular $p+q$ cube represented by the map 
$f \in {\mathcal{F}}_{\mathbf{S}}(U)$. Then $f^{\ast }(\alpha ) \in 
F^p(U)$ and $f^{\ast }(\beta )  \in F^q(U)$. Moreover, 
\begin{align}
\Lambda (\alpha , \beta ) \sigma & = \int_{J^{p+q}} f^{\ast }(\alpha ) \wedge f^{\ast }(\beta ) . 
\tag*{$\square $}
\end{align} 
We write $\alpha \wedge \beta $ instead of $\Lambda (\alpha , \beta )$. \medskip 

Next we define the \emph{exterior derivative} operator $\dee\, $. There is a unique real linear mapping 
\begin{equation}
\dee : G^p(\mathbf{S}) \rightarrow G^{p+1}(\mathbf{S}): 
{\lambda }_p(f_0, \ldots , f_p) \mapsto {\lambda }_{p+1}(\mathbf{1}, f_0, \ldots , f_p), 
\label{eq-s4.2one}
\end{equation}
where $f_i \in G^0(\mathbf{X})$. \medskip 

\noindent \textbf{Proposition 5.2.2} The map $\dee $ (\ref{eq-s4.2one}) is the coboundary operator with 
${\dee }^{\, 2} =0$. \medskip 

\noindent \textbf{Proof.} For $f_i \in G^0(\mathbf{S}) = {\mathcal{F}}_{\mathbf{S}}$ we have 
\begin{displaymath}
{\dee }^{\, 2} {\lambda }_p(f_0, \ldots , f_p) =\dee {\lambda }_{p+1}(\mathbf{1}, f_0, \ldots  ,f_p) = 
{\lambda }_{p+2}(\mathbf{1}, \mathbf{1}, f_0, \ldots , f_p) .
\end{displaymath}
From equation (\ref{eq-s4twonw}) we get 
\begin{displaymath}
{\lambda }_{p+2}(\mathbf{1}, \mathbf{1}, f_0, \ldots , f_p) = 
\int_{J^{p+2}} \mathbf{1} \frac{ \partial (\mathbf{1}, g_0, \ldots , g_p)}{\partial (t_1, \ldots , t_{p+1})} \, 
\dee t_1 \cdots \dee t_{p+1} , 
\end{displaymath}
where $g_i = f_i \comp f$ and $f \in {\mathcal{F}}_{\mathbf{S}}$ represents the singular 
$p+1$-cube $\sigma : J^{p+1} \rightarrow S$. Since $\frac{\dee \mathbf{1}}{\dee t_i} = 0$, 
the Jacobian $\frac{ \partial (\mathbf{1}, g_0, \ldots , g_p)}{\partial (t_1, \ldots , t_{p+1})} =0$. 
So ${\lambda }_{p+2}(\mathbf{1}, \mathbf{1}, f_0, \ldots , f_p) =0$. Consequently, 
${\dee }^{\, 2} =0$ on $G^p(\mathbf{S})$ for every $p \ge 0$. \hfill $\square $ \medskip 

From the definition of exterior product and exterior derivative we get  
\begin{equation}
{\lambda }_p(f_0, \ldots , f_p) = {\lambda }_0(f_0) \wedge {\lambda }_1(\mathbf{1}, f_1) \wedge \cdots 
{\lambda }_1(\mathbf{1}, f_p) = f_0 \, \dee f_1 \wedge \cdots \wedge \dee f_p , 
\label{eq-s4.2two}
\end{equation}
since ${\lambda }_0(f_0) = f_0$. It is straightforward to verify that 
\begin{displaymath}
\dee \, (\alpha \wedge \beta ) = \dee \alpha \wedge \beta + (-1)^{\mathrm{deg}\, \alpha } \alpha \wedge \dee \beta .
\end{displaymath}

Thus $\big( G(\mathbf{S}) = \sum_{p \ge 0} \oplus G^p(\mathrm{\mathbf{S}}), \dee \big) $ is a (graded) \emph{differential exterior algebra}. If $F: \mathbf{S} \rightarrow {\mathbf{S}}'$ is a differentiable map 
in the local category $\mathcal{D}$ of differential spaces and maps, then 
\begin{displaymath}
F^{\ast }: G({\mathbf{S}}') \rightarrow G(\mathbf{S}): \alpha \mapsto F^{\ast}(\alpha )
\end{displaymath} 
is a homomorphism of (graded) differential exterior algebras. Let 
$\mathcal{A}$ be the category of graded differential algebras and homomorphisms. Then the mapping $G: \mathcal{D} \rightarrow \mathcal{A}$, which sends the differentiable space 
$\mathbf{S}$ to the exterior algebra $G(\mathbf{S})$ and the mapping 
$F: {\mathbf{S}} \rightarrow {\mathbf{S}}'$ to the homomorphism 
$F^{\ast }: G({\mathbf{S}}') \rightarrow G(\mathbf{S})$ is a contravariant functor. \medskip 

When $\mathbf{S} = {\mathbf{S}}' \times {\mathbf{S}}''$ the space 
${\overline{\mathcal{F}}}_{{\mathbf{S}}' \times {\mathbf{S}}''}$, see page \pageref{pg-ten}, generates an exterior differential subalgebra 
$\overline{G}({\mathbf{S}} \times {\mathbf{S}}'')$ of $G({\mathbf{S}}' \times {\mathbf{S}}'')$. If 
the differentiable map $h: {\mathbf{S}}' \times {\mathbf{S}}'' \rightarrow \mathbf{T}$ is well behaved, then it induces a real linear mapping $h^{\ast }: G(\mathbf{T}) \rightarrow 
\overline{G}({\mathbf{S}}' \times {\mathbf{S}}'')$. \medskip 

Finally we prove the \emph{chain rule}. For $i=1, \ldots ,n$ let $g_i \in {\mathcal{F}}_{\mathbf{S}}$ 
and let  $F \in C^{\infty}({\R }^n)$. For $\phi = F(g_1, \ldots , g_n) \in {\mathcal{F}}_{\mathbf{S}}$ we have 
\begin{equation}%
\dee \phi  = \sum^n_{i=1} \frac{\partial F}{\partial g_i} \, \dee g_i. 
\label{eq-s4.2three}
\end{equation}

\noindent \textbf{Proof.} Let $\sigma : J^1 \rightarrow S$ be a singular $1$-cube represented 
by $f \in {\mathcal{F}}_{\mathbf{S}}(U)$. Then 
\begin{align}
\dee \phi (\sigma ) & = {\lambda }_1(\mathbf{1}, \phi ) \sigma = 
\int_{J^1} \frac{\dee \, (\phi  \comp f)}{\dee t} \, \dee t 
\notag \\
& = \int_{J^1} \frac{\dee }{\dee t} F (g_1 \comp f, \ldots , g_n \comp f) \, \dee t 
= \int_{J^1} \sum^n_{i=1} \frac{\partial F}{\partial g_i} \, \frac{\dee \, (g_i \comp f)}{\dee t} \, 
\dee t \notag \\
& = \sum^n_{i=1} \frac{\partial F }{\partial g_i} \, \int_{J^1} \frac{\dee \, (g_i \comp f)}{\dee t} \, \dee t = 
\sum^n_{i=1} \frac{\partial F }{\partial g_i} \, {\lambda }_1(\mathbf{1}, g_i) \sigma \notag \\
& = \big( \sum^n_{i=1} \frac{\partial F}{\partial g_i} \, \dee g_i \big) \sigma . \tag*{$\square $} 
\end{align} 

\subsection{The Poincar\'{e} lemma}

The unit interval $I = [0,1] \subseteq \R $ with $C^{\infty}(I) = C^{\infty}(\R )|_I$ is a differentiable 
space $\mathbf{I} = (I, C^{\infty}(I))$. A differentiable space $\mathbf{S} = 
(S, {\mathcal{F}}_{\mathbf{S}})$ 
is \emph{differentiably contractible} if there is a well behaved differentiable map 
$h: \mathbf{I} \times \mathbf{S} \rightarrow \mathbf{S}$ such that $h(1,x) =x $ and $h(0,x) = x_0$ for all $x \in S$ with $x_0$ a point in $S$. The map $h$ is called a \emph{contraction} of $S$. \medskip 

\noindent \textbf{Proposition 5.3.1} Let $\mathbf{S}$ be a differentiably contractible differentiable space and let $\alpha \in G^{p+1}(\mathbf{S})$ such that $\dee \alpha =0$. Then $\alpha = 
\dee \beta $ for some $\beta \in G^p(\mathbf{S})$. \medskip 

\noindent \textbf{Proof.} For $p \ge 0$ let $K_p(\mathbf{S})$, the set of \emph{singular $p$-chains} of $\mathbf{S}$, be the free abelian group generated by $C_p(\mathbf{S})$, the set of singular $p$-cubes of $\mathbf{S}$. Here $K_0(\mathbf{S}) = S$. Let $\sigma : J^{p+1} \rightarrow S$ be an element of $C_{p+1}(\mathbf{S})$, where 
\begin{displaymath}
J^{p+1}  = \{ (t'_1, \ldots , t'_{p+1}) \in {\R }^{p+1} \setrule \, A^{-}_i \le t'_i \le A^{+}_i, \, \, 1 \le i \le p+1 \} 
\end{displaymath} 
and $A^{\pm }_i \in \R $. For $1 \le i \le p$ let
\begin{align}
J^{p+1}_i & = \{ (t_1, \ldots , t_p) \in {\R }^p \setrule  \notag \\ 
& \hspace{.5in} A^{-}_j \le t_j \le A^{+}_j, \, \, 1 \le j < i \, \, \& \, \, A^{-}_{j+1} \le t_j \le A^{+}_{j+1} \, \, i \le j \le p \} . \notag 
\end{align} 
Then $J^{p+1}_i$ is a $p$-cube. For $1 \le i \le p$ let 
\begin{displaymath}
{\phi }^{\pm}_i: {\R }^p  \rightarrow {\R }^{p+1}: (t_1, \ldots , t_p) \mapsto 
(t'_1, \ldots , t'_{p+1}) = (t_1, \ldots ,t_{i-1}, A^{\pm}_i, 
t_i, \ldots , t_p) . 
\end{displaymath}
So 
\begin{align}
{\phi }^{\pm }_i(J^{p+1}_i) & = \{ (t'_1, \ldots , t'_{p+1}) = (t_1, \ldots , t_{i-1}, A^{\pm }_i, t_i, \ldots , t_p) 
\in {\R }^{p+1} \setrule \notag \\
& \hspace{-.25in} \, A^{-}_j \le t'_j \le A^{+}_j , \, 1 \le t'_j \le i-1 
\, \, \& \, \, A^{-}_{j+1} \le t'_{j+1} \le A^{+}_{j+1} \, i \le j \le p \} , \notag 
\end{align}
is the face $\{ t'_i = A^{\pm }_i \} $ of $J^{p+1}$. The \emph{geometric boundary} of the 
$p+1$-cube $J^{p+1}$ is 
\begin{displaymath}
\partial J^{p+1} = \bigcup_{1 \le i \le p} (-1)^{i+1}[ {\phi }^{+}_i(J^{p+1}_i) - {\phi }^{-}_i(J^{p+1}_i)] , 
\end{displaymath}
which is the union of appropriately signed faces of $J^{p+1}$. \medskip 

Let $\sigma : J^{p+1} \rightarrow S$ be a singular $p+1$-cube, that is, $\sigma \in 
K_{p+1}(\mathbf{S})$. Define the \emph{boundary} operator $\partial : K_{p+1}(\mathbf{S}) \rightarrow K_p(\mathbf{S})$ by 
\begin{displaymath}
\partial \sigma  = 
\sum^p_{i=1} (-1)^{i+1}[ \sigma \comp {\phi }^{+}_i - \sigma \comp {\phi }^{-}_i].
\end{displaymath}
For every $(t_1, \ldots , t_p) \in {\R }^p$ we have  
\begin{align}
\partial \sigma (t_1, \ldots , t_p) & = 
\sum^p_{i=1} (-1)^{i+1} [ \sigma (t_1, \ldots , t_{i-1}, A^{+}_i, t_i, \ldots , t_p) \notag \\
& \hspace{.5in} - \sigma (t_1, \ldots , t_{i-1}, A^{-}_i, t_i, \ldots , t_p)] . \notag 
\end{align}
So 
\begin{displaymath}
\partial \sigma |_{\partial J^{p+1}} = 
\sum^p_{i=1} (-1)^{i+1} [ \sigma \comp {\phi }^{+}_i(J^{p+1}_i)- \sigma \comp 
{\phi }^{-}_i(J^{p+1}_i)] . 
\end{displaymath}
Since 
\begin{displaymath}
\partial J^{p+1} =\bigcup_{1 \le i \le p} \big( (-1)^{i+1} {\phi }^{+}_i(J^{p+1}_i) \cup 
(-1)^i {\phi }^{-}_i(J^{p+1}_i) \big) , 
\end{displaymath}
the map $\partial \sigma |_{\partial J^{p+1}}: \partial J^{p+1} \rightarrow S$ is the sum of the signed singular 
$p$-cubes \linebreak 
${\sigma }^{\pm}_i: {\phi }^{\pm}_i(J^{p+1}_i) 
\rightarrow S$ for $1 \le i \le p$, in $C_p(\mathbf{S})$. Hence $\partial \sigma : \partial J^{p+1} 
\rightarrow S$ is a \linebreak 
singular $p$-chain in $K_p(\mathbf{S})$. Thus the boundary operator $\partial $ is 
well defined. When $p=0$, then $\sigma \comp {\phi }^{\pm}_i$ are the image under $\sigma $ of the end points of $J^1$. The boundary operator $\partial $ is a group homomorphism because it maps a generator of the group $K_{p+1}(\mathbf{S})$ into an element of $K_p(\mathbf{S})$. \medskip 

\noindent \textbf{Lemma 5.3.2} ${\partial }^{\, 2} =0$. \medskip 

\noindent \textbf{Proof.} Using the definition of the boundary operator $\partial $, we get 
\begin{align}
\partial (\partial \sigma )(t_1, \ldots , t_p) & = 
\sum^p_{j=1} (-1)^{j+1} \Big[ \partial \sigma (t_1, \ldots , t_{j-1}, A^{+}_j, t_j, \ldots , t_p) \notag \\
&\hspace{.5in} - \partial \sigma (t_1, \ldots , t_{j-1}, A^{-}_j, t_j, \ldots , t_p) \Big] \notag 
\end{align}
\begin{align}
\hspace{1.25in}&\hspace{-1in} = \sum^p_{j=1}\Big[ \sum_{a = \{ +, - \}}  \sum^j_{i=1}a (-1)^{i+j} \big[ 
\sigma (t_1, \ldots , t_{i-1}, A^{|a|}_i, t_i, \ldots , t_{j-1}, A^{a}_j, t_j, \ldots , t_p) \notag \\
& - \sigma (t_1, \ldots , t_{i-1}, A^{-|a|}_i, t_{i+1}, \ldots , t_{j-1}, A^{a}_j, t_j, \ldots , t_p) \big] \notag \\
& \hspace{-.5in} + \sum^p_{i =j} (-1)^i \big[ \sigma (t_1, \ldots , t_{j-1}, A^{a}_j, t_j, \ldots , t_{i-1}, 
A^{|a|}_i, t_i, \ldots , t_p) \notag \\
&- \sigma (t_1, \ldots , t_{j-1}, A^{a}_j, t_j, \ldots , t_{i-1}, A^{-|a|}_i, t_i, \ldots , t_p) \big]  \Big] , \notag 
\end{align} 
where $|a| = +$ if $a = +$ or $-$. 
The terms in the above expression, where the superscripts of $A^{|a|}_i$ and $A^a_j$ are the same sign,  cancel. We check this for $A^{+}_i$ and $A^{+}_j$.  We look at 
\begin{align}
&\sum^p_{j=1}\sum^j_{i=1}(-1)^{i+j}  
\sigma (t_1, \ldots , t_{i-1}, A^{+}_i, t_i, \ldots , t_{j-1}, A^{+}_j, t_j, \ldots , t_p) \notag \\
&- \sum^p_{j=1}\sum^p_{i =j} (-1)^{i+j} \sigma (t_1, \ldots , t_{j-1}, A^{+}_j, t_j, \ldots , t_{i-1}, A^{+}_i, t_i, \ldots , t_p) 
\label{eq-s4.3zero}
\end{align} 
Interchanging the indices, the second sum in (\ref{eq-s4.3zero}) reads 
\begin{displaymath}
- \sum^p_{j=1}\sum^j_{i =1} (-1)^{i+j} \sigma (t_1, \ldots , t_{i-1}, A^{+}_i, t_i, \ldots , t_{j-1}, A^{+}_j, t_j, \ldots , t_p) , 
\end{displaymath}
which cancels the first sum in (\ref{eq-s4.3zero}). \hfill $\square $ \medskip 

Let $\sigma : J^p \rightarrow S$ be an element of $C_p(\mathbf{S})$ with $p>0$. Define 
\begin{displaymath}
K(\sigma ): I \times J^p \rightarrow I \times S: (t'_1, t'_2, \ldots , t'_{p+1}) \mapsto 
(t'_1, \sigma (t'_2, \ldots , t'_{p+1})). 
\end{displaymath}
Then for $p >0$ we obtain 
\begin{displaymath}
K: K_p(\mathbf{S}) \rightarrow K_{p+1}(\mathbf{I} \times \mathbf{S}): \sigma \mapsto K(\sigma). 
\end{displaymath}
When $p=0$ we have $K(\sigma )(t_1) = (t_1, \sigma )$, where $\sigma \in K_0(\mathbf{S}) = S$. 
For $i=0,1$ define $u_i: S \rightarrow I \times S: x \mapsto (i,x)$. \medskip 

\noindent \textbf{Lemma 5.3.3}. We have 
\begin{align}
K \partial + \partial K & = (u_1)_{\ast } - (u_0)_{\ast } , 
\tag*{$(19)_p$}
\end{align}
where $(u_i)_{\ast }: K_p(\mathbf{S}) \rightarrow K_{p+1}(\mathbf{I} \times \mathbf{S})$ is the chain map induced by $u_i$. When $p=0$ equation $(19)_p$ reads 
\begin{align}
K \partial &= (u_1)_{\ast } - (u_0)_{\ast} . 
\tag*{$(19)_0$}
\end{align}

\noindent \textbf{Proof.} Let $\sigma : J^p \rightarrow S$ be a singular $p$-cube 
represented by $f \in \mathcal{F}(\mathbf{S})$. Then 
\begin{align}
K\partial \sigma (t_1, \ldots , t_p) + \partial K \sigma (t_1, \ldots , t_p) & = \notag \\
&\hspace{-1.85in} = K\big( \sum^p_{i=1}(-1)^{i+1}[ \sigma (t_1, \ldots , t_{i-1}, A^{+}_i, t_i , \ldots , t_p) 
- \sigma (t_1, \ldots , t_{i-1}, A^{-}_i, t_i , \ldots , t_p) ] \notag \\
& \hspace{-1in} + \partial (t_0, \sigma (t_1, \ldots , t_p) ) \big) \notag 
\end{align}
\begin{align}
\hspace{2in}& \hspace{-1.85in} = \sum^p_{i=1}(-1)^{i+1}[ \big( t_0, \sigma (t_1, \ldots , t_{i-1}, A^{+}_i, t_i , \ldots , t_p) \big) \notag \\
&\hspace{-1in} - \big( t_0, \sigma (t_1, \ldots , t_{i-1}, A^{-}_i, t_i , \ldots , t_p) \big) ] + 
\big( 1, \sigma (t_1, \ldots , t_p) \big) \notag \\
& \hspace{-1.5in} - \sum^p_{i=1} (-1)^{i+1}[ \big( t_0, \sigma (t_1, \ldots , t_{i-1}, A^{+}_i, t_i , \ldots , t_p) \big) \notag \\
&\hspace{-1in} - \big( t_0, \sigma (t_1, \ldots , t_{i-1}, A^{-}_i, t_i , \ldots , t_p) \big) ]  
- \big( 0, \sigma (t_1, \ldots , t_p) \big) \notag \\
& \hspace{-1.85in} = \big( 1, \sigma (t_1, \ldots , t_p) \big) - \big( 0, \sigma (t_1, \ldots , t_p) \big) \notag \\
& \hspace{-1.85in} = (u_1)_{\ast}\sigma (t_1, \ldots ,t_p) - (u_0)_{\ast }\sigma (t_1, \ldots , t_p). 
\tag*{$\square $} 
\end{align}

The set $K_p(\mathbf{S})$ of singular $p$-chains on $\mathbf{S}$, being the free abelian group 
generated by $C_p(\mathbf{S})$, the set of singular $p$-cubes on $\mathbf{S}$, may 
be considered as the set of vectors of a real vector space. The set $K^p(\mathbf{S})$ of singular 
$p$-cochains on $\mathbf{S}$ is the dual ${K_p(\mathbf{S})}^{\ast }$ of the vector space $K_p(\mathbf{S})$, because any function $f: K_p(\mathbf{S}) \rightarrow \R $ on the generators in $K_p(\mathbf{S})$ extends to a unique real linear function on the real vector space 
$K_p(\mathbf{S})$.  \medskip 

\noindent \textbf{Lemma 5.3.4} The dual operator 
\begin{displaymath}
K^{\ast }: K^{p+1}(\mathbf{S}) \rightarrow K^p(\mathbf{S}): \alpha \mapsto \alpha \comp K. 
\end{displaymath}
sends ${\overline{G}}^{p+1}(\mathbf{I} \times \mathbf{S})$ to $G^p(\mathbf{S})$. \medskip 

\noindent \textbf{Proof.} Every $f \in {\mathcal{F}}_{\mathbf{I} \times \mathbf{S}}$ can be 
represented as $f = \phi (t, f_1, \ldots , f_p)$, where $t \in I$ and $f_i \in {\mathcal{F}}_{\mathbf{S}}$. 
Every element $\omega $ of  ${\overline{G}}^{\, p+1}(\mathbf{I} \times \mathbf{S})$ is a sum of 
terms of the form 1: $a \dee t \wedge \dee f_1 \wedge \cdots \wedge \dee f_p$ or 
2: $a \dee f_1 \wedge \cdots \wedge \dee f_{p+1}$, where $a \in 
{\overline{\mathcal{F}}}_{\mathbf{I} \times \mathbf{S}}$ and $f_j \in {\mathcal{F}}_{\mathbf{S}}$. 
Let $\sigma : I \times J^p \rightarrow X$ be a singular $p+1$-cube in $S$ represented by $g: U  \rightarrow S$ in ${\mathcal{F}}_{\mathbf{S}}(U)$. Set $b = a \comp g$ and $g_j = f_j \comp g$. Suppose that $\omega $ has the form in item 2. Then $K^{\ast }\omega =0$. This follows because each $g_j$ is independent of $t$ and thus 
$\frac{\partial (g_1, \ldots , g_{p+1})}{\partial (t, t_1, \ldots , t_p)} =0$. When $\omega $ has the form in item 1, then 
\begin{align}
K^{\ast }\omega (\sigma ) & = \int_{I\times J^p} b \frac{\partial (g_1, \ldots , g_p)}{\partial (t_1, \ldots , t_p)} \, \dee t \dee t_1 \cdots \dee t_p \notag \\
& = \int_I b \, \dee t \, \int_{J^p} \frac{\partial (g_1, \ldots , g_p)}{\partial (t_1, \ldots , t_p)} \, \dee t_1 \cdots \dee t_{p+1}. \notag 
\end{align}
Consider $\int_I b \dee t$ as a function $f_0 : S \rightarrow \R $. Since $a \in 
{\overline{\mathcal{F}}}_{\mathbf{I} \times \mathbf{S}}$, it follows that $f_0 = \psi (f'_1, \ldots , f'_p)$ 
for some $f'_j \in {\mathcal{F}}_{\mathbf{S}}$ and some $\psi \in C^{\infty}({\R }^p)$. Thus 
$f_0 \in {\mathcal{F}}_{\mathbf{S}}$. Since 
$ \frac{\partial (g_1, \ldots , g_p)}{\partial (t_1, \ldots , t_p)}$ does not depend on $t$, it follows that 
$K^{\ast }\omega = f_0 \dee f_1 \wedge \cdots \wedge \dee f_p$. \hfill $\square $ \medskip  

\noindent \textbf{Lemma 5.3.5} For every $\mathbf{S} \in \mathcal{D}$, the dual 
${\partial }^{\ast }$ of the boundary operator $\partial :K_{p+1}(\mathbf{S}) \rightarrow 
K_p(\mathbf{S})$ is the exterior derivative $\dee : G^p(\mathbf{S}) \rightarrow 
G^{p+1}(\mathbf{S})$. \medskip 

\noindent \textbf{Proof.} Let $\sigma : J^{p+1} \rightarrow S$ be a singular $p+1$ cube represented by $f \in {\mathcal{F}}_{\mathbf{S}}(U)$. By definition, 
for $\alpha \in {K_p(\mathbf{S})}^{\ast } = K^p(\mathbf{S})$ we have 
\begin{align}
\alpha (\partial \sigma ) & = 
{\lambda }_p(f_0, \ldots , f_p)\partial \sigma , 
\quad \mbox{where $f_i \in {\mathcal{F}}_{\mathbf{S}}$} \notag \\
& = {\lambda }_p(f_0, \ldots , f_p) \big( \sum^p_{i=1}(-1)^{i+1} [ \sigma \comp 
{\phi }^{+}_i(J^{p+1}_i)- 
\sigma \comp {\phi }^{-}_i (J^{p+1}_i) \big) , \notag \\
&\hspace{.55in}\parbox[t]{2.25in}{since $\partial \sigma : \partial J^{p+1} \rightarrow S$ is a $p$-chain} \notag \\
& = \sum^p_{i=1}(-1)^{i+1}  \int_{{\phi }^{+}_i(J^{p+1}_i) - {\phi }^{-}_i(J^{p+1}_i)} 
g_0 \dee g_1 \wedge \cdots \wedge \dee g_p,  \notag  \\
& \hspace{.5in}\parbox[t]{2.25in}{where ${\lambda }_p(f_0, \ldots , f_p)
\big( \sigma \comp {\phi }^{\pm }_i(J^{p+1}_i) \big) =  $ \\ 
$ = \int_{{\phi }^{\pm }_i(J^{p+1}_i)} g_0 \dee g_1 \wedge \cdots \wedge \dee g_p$ and 
$g_i = $ $= f_i \comp f$, since $\partial \sigma $ is represented by $f$} \notag \\ 
& = \int_{\sum^p_{i=1}(-1)^{i+1}[{\phi }^{+}_i(J^{p+1}_i) - {\phi }^{-}_i(J^{p+1}_i)]} 
g_0 \dee g_1 \wedge \cdots \wedge \dee g_p \notag \\
& = \int_{\partial J^{p+1}} g_0 \dee g_1 \wedge \cdots \wedge \dee g_p \notag \\
& = \int_{J^{p+1}} \dee \big( g_0 \wedge \dee g_1 \wedge \cdots \wedge \dee g_p \big) , \, \, \mbox{by Stokes' theorem} \notag \\
& = {\lambda }_{p+1}(\mathbf{1}, f_0, \ldots , f_p) \sigma = \dee {\lambda }_p(f_0, \ldots , f_p)\sigma = \dee \alpha (\sigma) \notag  
\end{align} 
This implies ${\partial }^{\ast } = \dee $ because 
$({\partial }^{\ast }\alpha ) \sigma = \alpha (\partial \sigma ) = \dee \alpha (\sigma )$ 
for every $\sigma \in C^p(\mathbf{S})$. \hfill $\square $ \medskip 

We now finish the proof of the Poincar\'{e} lemma proposition 5.3.1. \medskip 

The real vector subspace ${\overline{G}}(\mathbf{I} \times \mathbf{S})$ of 
$G(\mathbf{I} \times \mathbf{S})$ 
is invariant under $\dee $. Hence the operator $(K \partial + \partial K)^{\ast } = \dee K^{\ast} + K^{\ast }\dee $ sends ${\overline{G}}^{\, p+1}(\mathbf{I} \times \mathbf{S})$ to $G^{p+1}(\mathbf{S})$. From 
equation $(19)_p$ we obtain 
\begin{align}
\dee K^{\ast } + K^{\ast }\dee & = u^{\ast }_1 - u^{\ast }_0 
\tag*{$(20)_p$}
\end{align}
on ${\overline{G}}^{\, p+1}(\mathbf{S})$. When $p=0$ equation $(20)_p$ reads 
\begin{align}
K^{\ast }\dee = u^{\ast }_1 - u^{\ast}_0. 
\tag*{$(20)_0$}
\end{align}
Let $h: I \times S \rightarrow S$ be a differentiable contraction of $S$. 
For $\alpha \in G^{p+1}(\mathbf{S})$ with $\dee \alpha =0$ set $\omega = h^{\ast }\alpha $. 
Since $h$ is a well behaved differentiable map, it follows that $\omega \in 
{\overline{G}}^{\, p+1}(\mathbf{I} \times \mathbf{S})$. Moreover, $\dee \omega = \dee h^{\ast }\alpha = h^{\ast }\dee \alpha =0$ because $\dee $ commutes with $h^{\ast }$. 
From $h \comp u_1 = {\mathrm{id}}_S$, we get $u^{\ast }_1 \omega = \alpha $. 
Since $h \comp u_0$ is a 
constant map, we obtain $u^{\ast }_0\omega =0$. Thus equation $(20)_p$ gives
\begin{displaymath}
\alpha = u^{\ast }_1 h^{\ast} \alpha - u^{\ast }_0 h^{\ast }\alpha = \dee \, (K^{\ast }\omega ) + 
K^{\ast }\dee \omega = \dee \, (K^{\ast }\omega ), 
\end{displaymath}
using $\dee \omega =0$. Setting $\beta = K^{\ast }\omega $ proves proposition 5.3.5. \hfill $\square $ \medskip 

\noindent \textbf{Corollary 5.3.5A} Suppose that the differentiable space $\mathbf{S}$ is 
differentiably contractible and $f \in {\mathcal{F}}_{\mathbf{S}}$ such that $\dee f =0$. 
Then $f$ is a constant function. \medskip 

\noindent \textbf{Proof.} This follows immediately from equation $(20)_0$ for 
\begin{displaymath}
u^{\ast }_1 f - u^{\ast }_0f = K^{\ast}(\dee f) = 0. 
\end{displaymath}
But $u^{\ast }_1 f = f$ and $u^{\ast }_0f$ is a constant function on $S$. \hfill $\square $ 

\subsection{Sheaf of differential forms}

Consider the differentiable space $\mathbf{S} = (S, \mathcal{F})$. Let $U$ be an open subset of 
$S$. Let $\sigma : J^p \rightarrow S$ be a singular $p$-cube on $\mathbf{S}$ represented by 
$f \in {\mathcal{F}}_{\mathbf{S}}$. Suppose that $\sigma (J^p) \subseteq U$. Then 
${\sigma }_U: J^p \rightarrow U: \mathbf{t} \mapsto \sigma (\mathbf{t})$ is a singular $p$ cube on $U$ represented by $f|_U \in {\mathcal{F}}_{U}$. Let $C_p(U)$ be the set of all 
singular $p$-cubes on $U$. Let $C^p(U)$ be the set of all functions 
$\alpha |_U : C_p(U) \rightarrow \R$ such that $\alpha : C_p(\mathbf{S}) \rightarrow \R $ is a singular $p$-cube on $\mathbf{S}$. For every $p > 0$ define the mapping 
\begin{displaymath}
({\lambda}_p)_U: {\mathcal{F}}^{p+1}_U \rightarrow C^p(U): 
(f_0|_U , \ldots , f_p|_U ) \rightarrow {{\lambda }_p}|_U (f_0|_U, \ldots , f_p|_U), 
\end{displaymath}
where $g_i = f|_U \comp f_i|_U \in {\mathcal{F}}_U$ and 
\begin{displaymath}
({\lambda }_p)_U (f_0|_U, \ldots , f_p|U)({\sigma }_U) = 
\int_{J^p} g_0 \frac{\partial (g_1, \ldots , g_p)}{\partial (t_1, \ldots , t_p)} \, \dee t_1 \cdots \dee t_p .
\end{displaymath}
When $p=0$ then $C^0(U) = {\mathcal{F}}|_U$. Let $G^0(U)$ be the vector space generated 
by the image of ${\lambda }_p|_U$. If $F: \mathbf{S} \rightarrow {\mathbf{S}}'$ is a differentiable 
maping such that $U =F^{-1}(U')$, where $U'$ is an open subset of ${\mathbf{S}}'$ and 
$U$ is an open subset of $\mathbf{S}$, then $F$ induces a real linear mapping 
\begin{displaymath}
F^{\ast }: C^p(U') \rightarrow C^p(U): ({\lambda }_p)_{U'}(f'_0, \ldots , f'_p) \mapsto 
({\lambda }_p)_U (f_0, \ldots , f_p), 
\end{displaymath}
where $f'_i = f_i \comp F \in {\mathcal{F}}_{U'}$ since $f_i \in {\mathcal{F}}_U$. Hence 
$F^{\ast }: G^p(U') \rightarrow G^p(U)$. Let ${\mathcal{G}}^p$ be the local category whose 
objects are the spaces $G^p(U)$ for $U$ an open subset of $\mathbf{S}$ and whose 
maps are real linear mappings $F^{\ast }: G^p(U') \rightarrow G^p(U)$ induced by 
differentiable mappings $F_U: U \rightarrow U'$. Thus  
$G: \mathcal{D} \rightarrow {\mathcal{G}}^p$ is a contravariant functor from the local category 
$\mathcal{D}$ of differentiable spaces and differentiable maps to the local category 
${\mathcal{G}}^p$ of differential $p$ forms on a differentiable 
space and real linear maps between spaces of differential $p$-forms. For every object $\mathbf{S}$
in the local category $\mathcal{D}$ we have a covariant functor $S_{\mathbf{S}}: 
{\mathcal{C}}^0 \rightarrow {\mathcal{G}}^p$ whose value at the object ${\mathbf{S}}_U$ is 
the open subset $U$ of $S$ and whose value at the inclusion map ${\iota }_U|V : V \subseteq U 
\rightarrow S$ is the map $({\iota }_U|_V)^{\ast }: G^p(V) \rightarrow G^p(U)$. For every 
$\mathbf{S}\in \mathcal{D}$ the map
\begin{displaymath}
G^p = G \comp S_{\mathbf{S}}: {\mathcal{C}}^0 \rightarrow {\mathcal{G}}^p: U \subseteq S \rightarrow G^p(U) 
\end{displaymath}
is the \emph{sheaf of differential $p$ forms} on $\mathbf{S}$. Note that the sheaf 
$G^0$ is equal to the sheaf $F^0$. 

\subsection{de Rham's theorem}

In this subsection we prove Smith's version of de Rham's theorem for  
differentiable spaces, see \cite{smith}. \medskip 

We now assume that the differentiable space $\mathbf{S}$ is \emph{locally differentiably contractible}, that is, for every $x \in S$ there is open subset $U$ of $S$ containing $x$ such that the differential space $(U, C^{\infty}(U))$ is contractible. We have the sheaf homomorphism 
$\dee : G^p \rightarrow G^{p+1}$ induced by the exterior derivative. Its 
kernel defines the sheaf $Z^p: {\mathcal{C}}^0 \rightarrow {\mathcal{G}}^p: U \mapsto 
G^p(U)$ of \emph{closed} $p$-forms on $S$.  \medskip 

\noindent \textbf{Lemma 5.5.1} The sheaf $Z^0$ is the sheaf $R$ of locally constant real valued 
functions on $\mathbf{S}$. \medskip 

\noindent \textbf{Proof.} Because $\mathbf{S}$ is locally differentiably contractible, from corollary 5.3.5A  it follows that if $\dee f =0$ on a contractible (and hence connected) open subset $U$ of 
$\mathbf{S}$, then $f$ is constant on $U$. \hfill $\square $ \medskip 

\noindent Hence we have an exact sequence of sheaves
\begin{equation}
0 \longrightarrow Z^q \stackrel{\iota }{\longrightarrow } {\Lambda }^q 
\stackrel{\dee }{\longrightarrow}Z^{q+1} \longrightarrow 0,  
\label{eq-4.5one}
\end{equation} 
for $q \ge 0$. Here ${\Lambda }^q:{\mathcal{C}}^0 \rightarrow {\mathcal{G}}^p: U \mapsto 
{\Lambda }^p(U)$ is the sheaf of differential $q$-forms on $\mathbf{S}$ 
generated by $f_i \in {\mathcal{F}}_{U}$, where $U$ is an open subset of $\mathbf{S}$.
In other words, $\alpha \in {\Lambda }^p(U)$ if and only if $\alpha = f_0 \dee f_1 \wedge \cdots 
\wedge \dee f_p$, where $f_i \in {\mathcal{F}}_U$. \medskip  

Let $\mathbf{S} \in \mathcal{D}$. Suppose that $S$ is a paracompact topological space. Then 
there is a locally finite \emph{good covering} $\mathcal{U} ={\{ U_i \} }_{i\in I}$ of $S$, that is, 
every finite intersection of elements of $\mathcal{U}$ is either contractible or empty. Consider 
the sheaf $F^0: {\mathcal{C}}^0 \rightarrow {\mathcal{F}}^0: U \mapsto {\mathcal{F}}_U$, 
where ${\mathcal{F}}^0$ is the category, for which an object is the space ${\mathcal{F}}_U$ of differentiable functions on the differentiable space ${\mathbf{S}}_U = (U, {\mathcal{F}}_U)$ with 
$U$ being an open subset of $S$,  and a map is 
\begin{displaymath}
({\iota }_U|_V)^{\ast}: {\mathcal{F}}_V \rightarrow {\mathcal{F}}_U: f \mapsto f \comp ({\iota }_U|_V), 
\end{displaymath}
We suppose that the sheaf $F^0$ is \emph{soft}, that is, 
for every closed subset $C$ of $S$ and every $g \in {\mathcal{F}}_{C} = 
{\mathcal{F}}_{\mathbf{S}}|_C$ there is an $f \in {\mathcal{F}}_{\mathbf{S}}$ such that $f|_C = g$. We have \medskip 

\noindent \textbf{Lemma 5.5.2} 
Let $\mathcal{S} = (S, {\mathcal{F}}_U)$ be a differentiable space. The sheaf ${\Lambda }^q$ of 
differential forms generated by ${\mathcal{F}}_U$ and exterior derivative of elements of 
${\mathcal{F}}_U$ for every open subset $U$ of $S$ is soft. \medskip

\noindent \textbf{Proof.} Let $\alpha = f_0 \dee f_1 \wedge \cdots \wedge \dee f_p \in 
{\Lambda }^p(C)$, where $f_i \in {\mathcal{F}}_C$ for some closed subset $C$ of $S$. 
Because the sheaf $F^0$ is soft, for each $f_i$ there is a $g_i \in {\mathcal{F}}_{\mathbf{S}}$ 
such that $g_i|_C = f_i$. Thus $\beta = g_0 \wedge \dee g_1 \wedge \cdots 
\wedge \dee g_p $ is a $p$-form on $S$, namely, $\beta \in {\Lambda }^p(S)$, such that 
$\beta |_C = \alpha $. Hence the sheaf ${\Lambda }^p$ is soft. \medskip 

We now review the definition of the $q^{\mathrm{th}}$ \emph{cohomology} group 
${\mathrm{H}}^q(\mathcal{U}, {\Lambda }^q)$ of the good covering $\mathcal{U}$ of 
$S$ having values in the sheaf ${\Lambda }^q$, see Lukina, et al. \cite{lukina-takens-broer}. 
Recall that a $q$-\emph{cochain} is an alternating map 
\begin{displaymath}
f: I^{q+1} \rightarrow {\Lambda }^q(U_{i_0} \cap \cdots \cap U_{i_q}) : (i_0, \ldots , i_q) 
\mapsto f(i_0, \ldots , i_q) ,
\end{displaymath}
where $U_{i_j} \in \mathcal{U}$ for $j = 0, \ldots , q$. The set $C^q(\mathcal{U}, {\Lambda }^q)$ of $q$-cochains is a group. The $q$-\emph{coboundary 
homomorphism} is the map 
\begin{displaymath}
\begin{array}{l}
{\delta }^q: C^q(\mathcal{U},{\Lambda }^q) \rightarrow C^{q+1}(\mathcal{U}, {\Lambda }^q): 
f \mapsto ({\delta }^qf)(i_0, \ldots , i_{q+1}) = \\
\hspace{.25in}= \sum^{q+1}_{j=0} {\iota }_{U_{i_0} \cap \cdots \cap \widehat{U}_{i_j} \cap \cdots \cap U_{i_{q+1} }} |_{U_{i_0} \cap \cdots \cap U_{i_{q+1}}}  \big( f( i_0, \ldots , i_{j-1}, \widehat{i_j}, i_{j+1}, \ldots , i_{q+1}) \big) .
\end{array}
\end{displaymath}
It is straightforward to show that ${\delta }^{q+1} \comp {\delta }^q =0$. We will omit the superscipt and 
write $\delta $ instead of ${\delta }^q$ when there is no risk of confusion. An element of $\ker {\delta }^q$ is called a $q$-\emph{cocycle}. The set of $q$-cocycles forms a subgroup of $C^q(\mathcal{U}, {\Lambda }^q)$. The group ${\mathrm{H}}^q(\mathcal{U}, {\Lambda }^q) = \ker {\delta }^q/ \mathrm{image}\, {\delta }^{q-1}$, when $q \ge 1$ and ${\mathrm{H}}^0(\mathcal{U}, {\Lambda }^q) = \ker {\delta }^0$, when $q =0$, 
is the $q^{\mathrm{th}}$ \emph{cohomology group} of the open covering $\mathcal{U}$ with 
values in the sheaf ${\Lambda }^q$. Note that ${\mathrm{H}}^0(\mathcal{U}, {\Lambda }^q) = 
Z^q(\mathbf{S})$. \medskip 

The following diagram of sheaf cochain complexes is exact. \medskip 

\hspace{-20pt} \begin{tabular}{ccccccccccc}
& & $\vdots $ & & $ \vdots $ & & $\vdots $ & & \\
& & $\downarrow $ & & $\downarrow $ & & $\downarrow $ & & \\
$0$ & $\longrightarrow $&$C^{k-1}(\mathcal{U}, Z^q) $& $\stackrel{\iota}{\longrightarrow}$ & 
$C^{k-1}(\mathcal{U}, {\Lambda }^q)$ &$\stackrel{\dee}{\longrightarrow}$ & $C^{k-1}(\mathcal{U},Z^{q+1})$ & 
$\longrightarrow $& $0$ \\
& & $\rule{5pt}{0pt} \downarrow $ $\delta $ & & $\rule{6pt}{0pt} \downarrow$ $\delta $ & & $\rule{5pt}{0pt} \downarrow$ $ \delta $ & & \\
$0$ & $\longrightarrow $&$C^k(\mathcal{U}, Z^q) $& $\stackrel{\iota}{\longrightarrow}$ & 
$C^k(\mathcal{U}, {\Lambda }^q)$ &$\stackrel{\dee}{\longrightarrow}$ & $C^k(\mathcal{U},Z^{q+1})$ & 
$\longrightarrow $& $0$ \\
& & $\downarrow $ & & $\downarrow $ & & $\downarrow $ & & \\
& & $\vdots $ & & $\vdots $ & & $\vdots $ & &
\end{tabular}
\vspace{.1in}
\par \noindent where $\delta $ is the coboundary operator, $\iota $ is the inclusion, and $\dee $ the exterior derivative. The above diagram of cochain complexes gives rise to the long exact sequence of cohomology groups \medskip 
\begin{equation}
0 \rightarrow {\mathrm{H}}^0(\mathcal{U}, Z^q) \stackrel{{\iota }_{\ast }}{\longrightarrow} 
{\mathrm{H}}^0(\mathcal{U}, {\Lambda }^q) \stackrel{{\dee }_{\, \ast}}{\longrightarrow} 
{\mathrm{H}}^0(\mathcal{U}, Z^{q+1}) \stackrel{{\delta }_{\ast }}{\longrightarrow} 
{\mathrm{H}}^1(\mathcal{U}, Z^q) \rightarrow \cdots , 
\label{eq-s5.2one}
\end{equation} 
where ${\iota }_{\ast }$, ${\dee }_{\, \ast }$, and ${\delta }_{\ast }$ are cohomology homomorphisms 
induced by $\iota $, $\dee $\, , and $\delta $ respectively. Because the sheaf ${\Lambda }^q$ is 
soft, its $q^{\mathrm{th}}$ cohomology group vanishes when $q \ge 1$. Thus the long exact 
sequence (\ref{eq-s5.2one}) falls apart into the exact sequences 
\begin{equation}
0 \rightarrow {\mathrm{H}}^0(\mathcal{U}, Z^q) \stackrel{{\iota }_{\ast}}{\longrightarrow} 
{\mathrm{H}}^0(\mathcal{U}, {\Lambda }^q) \stackrel{{\dee }_{\, \ast}}{\longrightarrow} 
{\mathrm{H}}^0(\mathcal{U}, Z^{q+1}) \stackrel{{\delta }_{\ast }}{\longrightarrow} 
{\mathrm{H}}^1(\mathcal{U}, Z^q) \rightarrow 0 
\label{eq-s4.5two}
\end{equation}
and 
\begin{equation}
0 \rightarrow {\mathrm{H}}^k(\mathcal{U}, Z^{q+1}) \stackrel{{\delta }_{\ast }}{\longrightarrow} 
{\mathrm{H}}^{k+1}(\mathcal{U}, Z^q) \rightarrow 0,  
\label{eq-s4.5three}
\end{equation}
for every $k \ge 1$. Now ${\mathrm{H}}^0(\mathcal{U}, R) = Z^0(\mathbf{S})$. Applying 
(\ref{eq-s4.5three}) iteratively gives
\begin{displaymath}
{\mathrm{H}}^q(\mathcal{U}, R) = {\mathrm{H}}^q(\mathcal{U}, Z^0) \simeq 
{\mathrm{H}}^{q-1}(\mathcal{U}, Z^1) \simeq \cdots \simeq {\mathrm{H}}^1(\mathcal{U}, Z^{q-1}), 
\end{displaymath}
where $\simeq $ means isomorphic to. For $q \ge 1$ from (\ref{eq-s4.5two}) we find that 
\begin{align}
{\mathrm{H}}^1(\mathcal{U}, Z^{q-1}) & \simeq 
{\mathrm{H}}^0(\mathcal{U}, Z^q)/\ker {\delta }_{\ast }|_{{\mathrm{H}}^0(\mathcal{U}, Z^q)} 
\notag \\
& \simeq 
{\mathrm{H}}^0(\mathcal{U}, Z^q)/ {\dee }_{\, \ast } \big( {\mathrm{H}}^0(\mathcal{U}, {\Lambda }^{q-1}) 
\big) \simeq Z^q(\mathbf{S})/\dee {\Lambda }^{q-1}(\mathbf{S}).  \notag
\end{align}
Hence we obtain Smith's version of de Rham's theorem for a differential space. \medskip

\noindent \textbf{Theorem 5.5.3} Let $\mathbf{S} = (S, {\mathcal{F}}_{\mathbf{S}})$ be 
a differentiable space, where $S$ is locally 
differentiably contractible and paracompact. Moreover, suppose that the sheaf 
$F^0:{\mathcal{C}}^0 \rightarrow {\mathcal{F}}^0: U \mapsto {\mathcal{F}}_U$ is soft. Then for 
every $q \ge 1$ \medskip 
\begin{align}
&{\mathrm{H}}^q(\mathcal{U}, R) = {\mathrm{H}}^q(\mathcal{U}, Z^0) \simeq  
Z^q(\mathbf{S})/\dee {\Lambda}^{q-1}(\mathbf{S}) = {\mathrm{H}}^q_{\mathrm{dR}}(\mathbf{S}) 
\tag*{$(23)_q$}
\end{align}
and for $q=0$ 
\begin{align}
& {\mathrm{H}}^0(\mathcal{U}, R) = Z^0(\mathbf{S}) = R = {\mathrm{H}}^0_{\mathrm{dR}}(\mathbf{S}).
\tag*{$(23)_0$}
\end{align} 

\noindent \textbf{Corollary 5.5.3A} For $q \in {\Z }_{\ge 0}$ the $q^{\mathrm{th}}$ 
\v{C}ech sheaf cohomology group ${\check{H}}^q (\mathbf{S}, R)$ of 
$\mathbf{S}$ with values in the sheaf $R$ of locally constant 
real valued functions on $\mathbf{S}$ is isomorphic to the $q^{\mathrm{th}}$ 
de Rham cohomology group ${\mathrm{H}}^q_{\mathrm{dR}}(\mathbf{S})$ of $\mathbf{S}$. \medskip

\noindent \textbf{Proof.} From equations $(23)_q$ and $(23)_0$ it follows that 
${\mathrm{H}}^q(\mathcal{U}, R) \simeq {\mathrm{H}}^q_{\mathrm{dR}}(\mathbf{S})$, 
which is independent of the good covering $\mathcal{U}$ of $\mathbf{S}$. By definition 
${\check{H}}^q(\mathbf{S}, R) = {\mathrm{H}}^q(\mathcal{U}, R)$. 
\hfill $\square $ 

\section{The orbit space of a proper action}

Let $\Phi: G \times M \rightarrow M:(g,m) \mapsto {\Phi }_g(m) = g \cdot m$ be a proper action of a 
Lie group $G$ on a smooth manifold $M$. We denote the space of 
smooth $G$-invariant functions on $M$ by ${C^{\infty}(M)}^G$. Let 
$M/G = \{ G \cdot m \setrule \, m \in M \} = \overline{M}$ be the space of $G$-orbits in $M$ with $G$-orbit map $\pi : M \rightarrow \overline{M}: m \mapsto \overline{m} = G\cdot m$. The set 
$C^{\infty}(\overline{M})= \{ f : \overline{M} \rightarrow \R \setrule {\pi }^{\ast }f \in {C^{\infty}(M)}^G \} $ 
is a differential structure on $\overline{M}$. The differential space 
$(\overline{M},C^{\infty}(\overline{M}))$ is subcartesian and is locally 
closed. Moreover, the differential space topology on 
$\overline{M}$ is equivalent to the quotient topology. So $\overline{M}$ is paracompact. 

\subsection{Vector fields on $\overline{M}$}

Let $X$ be a $G$-invariant vector field on $M$, that is, ${\Phi }^{\ast}_gX(m) =  
T_m{\Phi }_{g^{-1}}X({\Phi }_g(m))$ $= X(m)$ for every $g \in G$ and every $m \in M$. \medskip 

\noindent \textbf{Lemma 6.1.1} If $X$ is a $G$-invariant vector field on $M$, then 
\begin{displaymath}
X: C^{\infty}(M)^G \rightarrow C^{\infty}(M)^G: f \mapsto X(f) = L_Xf 
\end{displaymath}
is a derivation of $C^{\infty}(M)^G$. \medskip

\noindent \textbf{Proof.} $X$ maps $C^{\infty}(M)^G$ into itself, because 
for every $g \in G$ and every $f \in C^{\infty}(M)^G$ we have 
\begin{displaymath}
{\Phi }^{\ast }_g(L_Xf) = L_{{\Phi }^{\ast }_gX}\big( {\Phi }^{\ast}_g(f)\big) = L_Xf.
\end{displaymath}
$X$ is a derivation, because it is clearly a real linear mapping of $C^{\infty}(M)^G$ into 
itself, and it satisfies Leibniz' rule, since for $f_1$, $f_2 \in C^{\infty}(M)^G$ we have 
\begin{align}
X(f_1f_2)& = L_X(f_1f_2) = (L_Xf_1)f_2 +f_1(L_Xf_2) = X(f_1)f_2 + f_1X(f_2).
\tag*{$\square $}
\end{align}

\noindent For each $m \in M$ let ${\varphi }_m: I_m \subseteq \R \rightarrow M:t \mapsto 
{\varphi }_t(m)$ be a maximal integral curve of the $G$-invariant vector field $X$ on $M$, which starts at $m$. Here $I_m$ is an open interval containing $0$. Let 
\begin{displaymath}
\mathcal{D} = {\amalg}_{m\in M}I_m 
= \{ (t, m) \in \R \times M| \, t \in I_m \, \, \mbox{for each $m \in M$}  \} . 
\end{displaymath}
Then $\mathcal{D}$ 
is an open subset of $\R \times M$, which is the \emph{domain} of the \emph{local flow} 
$\varphi : \mathcal{D} \subseteq \R \times M \rightarrow M$ of the vector field $X$, namely, \medskip 

\indent 1) For each $m \in M$ we have ${\varphi }_0(m) = m$; \smallskip 

2) \parbox[t]{4.5in}{if $(t,m)\in \mathcal{D}$, $(s, {\varphi }_t(m))\in \mathcal{D}$, and 
$(s+t,m) \in \mathcal{D}$, then ${\varphi }_{s+t}(m) = {\varphi }_s\big( {\varphi }_t(m) \big) $.} \smallskip 

\vspace{-.1in}\noindent For each $m \in M$ there is an open neighborhood $U_m$ of $m$ in $M$ such that the map ${\mathcal{D}}_{U_m} \subseteq \R \rightarrow C^{\infty}(U_m, M): t \mapsto {\varphi }_t$, where ${\mathcal{D}}_{U_m} = \{ t \in \R \setrule \, (t,m') \in \mathcal{D} \, \, 
\mbox{for every $m' \in U_m$} \} $, is a local one parameter group of local 
diffeomorphisms of $M$. \medskip

\noindent \textbf{Lemma 6.1.2} $X$ is a $G$-invariant vector field on $M$ if and only if 
\begin{equation}
({\Phi }_g \comp {\varphi }_t)(m) = ({\varphi }_t \comp {\Phi }_g)(m)
\label{eq-s6ss1onevnw}
\end{equation}
for every $(t,m) \in \mathcal{D}$. \medskip 

\noindent \textbf{Proof.} If equation (\ref{eq-s6ss1onevnw}) holds, then differentiating 
both sides with respect to $t$ and then evaluating at $t =0$ gives the second equality below
\begin{align}
T_m{\Phi }_g X(m) & = T_m{\Phi }_g \dbydt \hspace{-8pt}{\varphi }_t (m), \quad 
\parbox[t]{2.5in}{since $t \mapsto {\varphi }_t(m)$ is an integral curve of $X$ 
starting at $m$.} \notag \\
& = \dbydt \hspace{-8pt}{\varphi }_t\big( {\Phi }_g(m) \big)  = X({\Phi }_g(m)). \notag 
\end{align}
Consequently, ${\Phi }^{\ast }_gX = X$, that is, the vector field $X$ is $G$-invariant. To 
prove the converse, note that $t \mapsto {\varphi }_t({\Phi }_g(m))$ is an integral 
curve of $X$ starting at ${\Phi }_g(m)$. We now show that $t \mapsto {\Phi }_g({\varphi }_t(m))$ 
is also an integral curve of $X$ starting at ${\Phi}_g(m)$. We have 
\begin{align}
\frac{\dee }{\dee t}({\Phi }_g \comp {\varphi }_t) (m) & = 
T_{{\varphi }_t(m)}{\Phi }_g \Big( \frac{\dee }{\dee t}{\varphi }_t(m) \Big) \notag \\
& = T_{{\varphi }_t(m)}{\Phi }_g \big( X({\varphi }_t(m)) \big) , \quad \parbox[t]{2in}{since 
$t \mapsto {\varphi }_t(m)$ is an integral curve of $X$ starting at $m$} \notag \\
& = X({\Phi }_g({\varphi }_t(m))), \quad \parbox[t]{2.5in}{since $X$ is a $G$-invariant vector 
field} \notag \\
& = X\big( ({\Phi }_g \comp {\varphi }_t)(m) \big). \notag 
\end{align}
Since maximal integral curves of the vector field $X$ starting at ${\Phi }_g(m)$ are unique, 
we get $({\Phi }_g \comp {\varphi }_t)(m)= ({\varphi }_t \comp {\Phi }_g)(m)$ for every 
$(t,m) \in \mathcal{D}$. \hfill $\square $ \medskip 

From equation (\ref{eq-s6ss1onevnw}) it follows that ${\varphi }_t(G \cdot m) = 
G\cdot {\varphi }_t(m)$ for every $(t,m) \in \mathcal{D}$. Hence ${\varphi }_t$ induces the map 
\begin{equation}
{\overline{\varphi }}_t: \overline{M} \rightarrow \overline{M}: 
\overline{m} = G \cdot m \mapsto G \cdot {\varphi }_t(m) = \overline{{\varphi }_t(m)}. 
\label{eq-s6ss1twovnw}
\end{equation}

\noindent \textbf{Lemma 6.1.3} For every $(t,m) \in \mathcal{D}$ we have 
\begin{equation}
{\overline{\varphi }}_t \comp \pi = \pi \comp {\varphi }_t.
\label{eq-s6ss1threevnw*}
\end{equation}

\noindent \textbf{Proof.} Equation (\ref{eq-s6ss1threevnw*}) follows from equation 
(\ref{eq-s6ss1twovnw}) because 
\begin{displaymath}
{\overline{\varphi }}_t\big( \pi (m) \big) = {\overline{\varphi }}_t(\overline{m}) 
=\overline{{\varphi }_t(m)} = \pi \big( {\varphi }_t(m) \big) , 
\end{displaymath}
for every $(t,m) \in \mathcal{D}$. \hfill $\square $ \medskip 

Let $X$ be a $G$-invariant vector field on $M$. Let 
\begin{equation}
\overline{X}: C^{\infty}(\overline{M}) \rightarrow C^{\infty}(\overline{M}): 
\overline{f} \mapsto \overline{X({\pi }^{\ast }\overline{f})}, 
\label{eq-s6ss1threevnw}
\end{equation}
that is, ${\pi}^{\ast }(\overline{X}(\overline{f})) = X({\pi }^{\ast }\overline{f})$, 
for every $\overline{f} \in C^{\infty}(M)$. We can rewrite the preceding equation as 
\begin{equation}
T\pi \comp X = \overline{X} \comp \pi . 
\label{eq-s6ss1fivevnw}
\end{equation}

\noindent \textbf{Lemma 6.1.4} $\overline{X}$ is a derivation of $C^{\infty}(\overline{M})$. \medskip 

\noindent \textbf{Proof.} $\overline{X}$ is a real linear mapping of $C^{\infty}(\overline{M})$ 
into itself, for if $\alpha $, $\beta \in \R $ and ${\overline{f}}_1$, ${\overline{f}}_2 \in 
C^{\infty}(\overline{M})$ then 
\begin{align}
\overline{X}(\alpha {\overline{f}}_1 + \beta {\overline{f}}_2) & = 
\overline{X({\pi }^{\ast}(\alpha {\overline{f}}_1 + \beta {\overline{f}}_2))}, 
\quad \mbox{by definition of $\overline{X}$} \notag \\
&\hspace{-.5in}  = \overline{X({\pi }^{\ast }(\alpha {\overline{f}}_1) + 
{\pi }^{\ast }(\beta {\overline{f}}_2))} \notag \\
&\hspace{-.5in} = \alpha \overline{X({\pi }^{\ast }({\overline{f}}_1))} + 
\beta \overline{X({\pi }^{\ast }({\overline{f}}_2))}, 
\, \,  \mbox{since $X$ is a derivation of $C^{\infty}(M)^G$} \notag \\
&\hspace{-.5in} = \alpha \overline{X}({\overline{f}}_1) + \beta \overline{X}({\overline{f}}_2). \notag
\end{align}
Also $\overline{X}$ satisfies Leibniz' rule, for 
\begin{align}
\overline{X}({\overline{f}}_1{\overline{f}}_2) & = 
\overline{X({\pi }^{\ast }({\overline{f}}_1 {\overline{f}}_2)} = 
\overline{X({\pi }^{\ast }{\overline{f}}_1\, {\pi }^{\ast }{\overline{f}}_2)} \notag \\
& = \overline{X({\pi }^{\ast }{\overline{f}}_1)} \, \overline{{\pi }^{\ast }{\overline{f}}_2 }+ 
\overline{{\pi }^{\ast }{\overline{f}}_1} \, \overline{X({\pi }^{\ast }{\overline{f}}_2)} \notag \\
& = \overline{X}({\overline{f}}_1){\overline{f}}_2 + {\overline{f}}_1 \overline{X}({\overline{f}}_2). 
\tag*{$\square $} 
\end{align}

\noindent \textbf{Proposition 6.1.5} Let $\overline{\mathcal{D}} = \{ (t, \overline{m}) \in \R \times \overline{M} \setrule \, 
t \in {\mathcal{D}}_m \, \, \& \, \, \overline{m} \in \pi ({\mathcal{D}}_t ) \} $ 
with ${\mathcal{D}}_t = \{ m \in M \setrule \, (t,m) \in \mathcal{D} \}$ and 
${\mathcal{D}}_m = \{ t \in \R \setrule \, (t, m) \in \mathcal{D} \}$. Then the 
curve
\begin{displaymath}
\overline{c}: {\overline{\mathcal{D}}}_{\overline{m}} = {\mathcal{D}}_m \subseteq \R \rightarrow \overline{M}: \overline{m} \mapsto {\overline{\varphi }}_t(\overline{m}) 
\end{displaymath}
is a maximal integral curve of the derivation $\overline{X}$ of $C^{\infty}(\overline{M})$ starting at 
$\overline{m} \in \overline{M}$. \medskip 

\noindent \textbf{Proof.} Let $c: {\mathcal{D}}_m \subseteq \R \rightarrow M : 
t \mapsto {\varphi }_t(m)$ be a maximal integral curve of the $G$-invariant vector field $X$ on $M$. Then 
\begin{equation}
Tc (t) = (X \comp c)(t)\quad \mbox{for every $t \in {\mathcal{D}}_m$ and every $m \in M$}. 
\label{eq-s6ss1sixvnw}
\end{equation}
Since $\overline{c} = {\pi }_{\ast }c$ we get 
\begin{align}
T\overline{c}(t) & = T({\pi }_{\ast} c)(t) = (T{\pi }_{\ast })(Tc) (t) \notag \\
& = (T{\pi }_{\ast })(X\comp c))(t), \quad \mbox{using (\ref{eq-s6ss1sixvnw})} \notag \\
& = \overline{X} ({\pi }_{\ast } c )(t), \quad \mbox{using (\ref{eq-s6ss1fivevnw})} \notag \\
& = (\overline{X} \comp \overline{c})(t). \notag
\end{align}
Thus $\overline{c}$ is an integral curve of the derivation $\overline{X}$, which is maximal 
because by equation (\ref{eq-s6ss1fivevnw}) the derivations $X$ on $C^{\infty}(M)^G$ and 
$\overline{X}$ on $C^{\infty}(\overline{M})$ are $\pi $-related and their local flows satisfy 
equation (\ref{eq-s6ss1threevnw}). \hfill $\square $ \medskip 

\noindent \textbf{Corollary 6.1.5A} The map 
\begin{displaymath}
\overline{\varphi }: \overline{\mathcal{D}} \subseteq \R \times \overline{M} \rightarrow \overline{M}:
(t, \overline{m}) \mapsto {\overline{\varphi}}_t(\overline{m})
\end{displaymath}
is the local flow of the derivation $\overline{X}$ on $C^{\infty}(\overline{M})$. \medskip 

\noindent \textbf{Proof.} This follows immediately from proposition 6.1.5. \hfill $\square $ \medskip 

\noindent \textbf{Corollary 6.1.5B} The derivation $\overline{X}$ on $C^{\infty}(\overline{M})$ is 
a vector field on $\overline{M}$. \medskip 

\noindent \textbf{Proof.} Since ${\mathcal{D}}_t$ is an open subset of $M$ for all $t \in \R $ and 
the $G$-orbit map $\pi $ is an open map, we get ${\overline{\mathcal{D}}}_t = 
\pi \big( {\mathcal{D}}_t \big) $ is an open 
subset of $\overline{M}$. Also ${\mathcal{D}}_m$ is a open interval containing $0$ for all 
$m \in M$, which implies ${\overline{\mathcal{D}}}_{\overline{m}}$ is an open interval 
containing $0$ for all $\overline{m} \in \overline{M}$, since ${\overline{\mathcal{D}}}_{\overline{m}} = {\mathcal{D}}_m$. Thus every maximal integral curve of $\overline{X}$ is defined on an 
open interval in $\R $ that contains $0$ and is also unique since every integral curve of $X$ is 
unique and the derivations $X$ and $\overline{X}$ are $\pi $-related. 
Hence by proposition 3.3.1 the derivation 
$\overline{X}$ on $(\overline{M}, C^{\infty}(\overline{M}))$ is a vector field. 
\hfill $\square $ 

\subsection{The differentiable space $\overline{M}$}

In this subsection we show that the $G$-orbit space $(\overline{M}, C^{\infty}(\overline{M}))$ is 
a differentiable space in the sense of Smith. \medskip 

\noindent \textbf{Lemma 6.2.1} The surjective continuous $G$-orbit mapping 
$\pi : M \rightarrow \overline{M}$ at every $\overline{m} \in \overline{M}$ has 
a continuous local section, that is, there is an open neighborhood $\overline{U}$ of 
$\overline{m}$ in $\overline{M}$ and a continuous mapping 
${\sigma }_{\overline{U}}: \overline{U} \subseteq \overline{M} \rightarrow M$ such that 
$\pi \comp {\sigma }_{\overline{U}} = {\mathrm{id}}_{\overline{U}}$. \medskip 

\vspace{-.15in}\noindent \textbf{Proof.} Let $\overline{m} \in \overline{M}$. 
Suppose that $\widetilde{U}$  is 
an open subset of $M$ containing $m \in {\pi }^{-1}(\overline{m})$. Since the 
$G$ action on $M$ is proper, it has a slice $S_m$ at $m$. Let 
$U$ be an open neighborhood of $m$, which is contained in $\widetilde{U}$, such that 
its closure $\mathrm{cl}(U)$ is compact. Then the mapping $\pi |_{\mathrm{cl}(U) \cap S_m}$ 
is a continuous bijection onto $\pi \big( \mathrm{cl}(U) \big) $, which is a compact 
subset of the Hausdorff topological space $\overline{M}$. Hence the map 
$\pi |_{\mathrm{cl}(U) \cap S_m}$ is a homeomorphism of $U \cap S_m$ onto the 
open neighborhood $\overline{U} = \pi (U \cap S_m)$ of $\overline{m} = \pi (m)$. 
Let ${\sigma }|_{\overline{U}} = ({\pi }^{-1})|_{\overline{U}}: \overline{U} \subseteq \overline{M} 
\rightarrow M $. Since $\pi \comp {\sigma }_{\overline{U}} = {\mathrm{id}}_{\overline{U}}$ by 
construction, ${\sigma }_{\overline{U}}$ is a continuous local section of the mapping 
$\pi $. \hfill $\square $ \medskip 

\noindent \textbf{Theorem 6.2.2} The differential space $(\overline{M}, C^{\infty}(\overline{M} ))$ is 
a differentiable space in the sense of Smith. \medskip 

To prove theorem 6.2.2 we will need some preliminary results. \medskip

\noindent \textbf{Lemma 6.2.3} If $\overline{c} \in C^0(\R , \overline{M})$, the set of 
continuous maps from $\R $ to the $G$-orbit space $\overline{M}$, then there is 
$c \in C^0(\R , M)$, which covers $\overline{c}$, that is, $\overline{c} = \pi \comp c$. \medskip 

\vspace{-.15in}\noindent \textbf{Proof.} Let $\overline{c} \in C^0(\R , \overline{M})$ be given. Let 
${\mathcal{T}}^{+}$ be the set of all $t \in {\R }_{> 0}$ such that there is a continuous 
curve $c:[0,t) \rightarrow M$, which covers the given continuous curve 
$\overline{c}|_{[0,t)}$. Let $c(0) \in {\pi }^{-1}(\overline{c}(0))$. By lemma 6.2.1 there is 
an open neighborhood ${\overline{U}}_{\overline{c}(0)}$ of $\overline{c}(0)$ in $\overline{M}$ 
and a local section ${\sigma }_{{\overline{U}}_{\overline{c}(0)}}: {\overline{U}}_{\overline{c}(0)} 
\subseteq \overline{M} \rightarrow M$ of the $G$-orbit map $\pi : M \rightarrow \overline{M}$ 
such that ${\sigma }_{{\overline{U}}_{\overline{c}(0)}}\big( \overline{c}(0) \big) = c(0)$. Since 
${\overline{U}}_{\overline{c}(0)}$ is an open neighborhood of $\overline{c}(0)$, there is 
an ${\varepsilon }_{\overline{c}(0)} > 0$ such that $\overline{c}\big( [0, {\varepsilon }_{\overline{c}(0)})\big) 
\subseteq {\overline{U}}_{\overline{c}(0)}$. On $[0, {\varepsilon }_{\overline{c}(0)})$ let 
$c = {\sigma }_{{\overline{U}}_{\overline{c}(0)}} \comp \overline{c}$. Then 
$c$ is a continuous curve starting at $c(0)$ and defined on $[0, {\varepsilon }_{\overline{c}(0)})$, 
which covers $\overline{c}$ on $[0, {\varepsilon }_{\overline{c}(0)})$, since 
$\pi \comp c = \pi \comp ({\sigma }_{{\overline{U}}_{\overline{c}(0)}} \comp c ) = c$ on 
$[0, {\varepsilon }_{\overline{c}(0)})$. This shows that ${\varepsilon }_{\overline{c}(0)} \in 
{\mathcal{T}}^{+}$. Since ${\mathcal{T}}^{+}$ is nonempty, it has a supremum $T^{+}$. 
Suppose that $T^{+} < \infty $. By lemma 6.2.1 there is an open neighborhood 
${\overline{U}}_{\overline{c}(T^{+})}$ of $\overline{c}(T^{+})$ in $\overline{M}$ and a 
continuous local section ${\sigma }_{{\overline{U}}_{\overline{c}(T^{+})}}: 
{\overline{U}}_{\overline{c}(T^{+})} \subseteq \overline{M} \rightarrow M$ of the $G$-orbit map 
$\pi $ such that for some $\widetilde{t} \in {\mathcal{T}}^{+}$ we have $\overline{c}(\widetilde{t}) 
\in {\overline{U}}_{\overline{c}(T^{+})}$. Since $\widetilde{t} \in {\mathcal{T}}^{+}$, there is 
a continuous curve $c$ on $M$ defined on $[0, \widetilde{t})$, which covers 
${\overline{c}}|_{[0, \widetilde{t})}$. So $\overline{c}(\widetilde{t}) = \pi (c(\widetilde{t}))$. 
Define $\widetilde{c} = {\sigma }_{{\overline{U}}_{\overline{c}(T^{+})}} \comp \overline{c}$. 
Then $\widetilde{c}$ is a continuous curve on $[0, T^{+} + {\varepsilon }_{\overline{c}(T^{+})}) $, 
where ${\varepsilon }_{\overline{c}(T^{+})} > 0$ is chosen so that 
$\overline{c}\big( [0, T^{+} + {\varepsilon }_{\overline{c}(T^{+})}) \big) \subseteq 
{\overline{U}}_{\overline{c}(T^{+})} $. Moreover, ${\widetilde{c}}|_{[0, \widetilde{t})} = 
c|_{[0, \widetilde{t})}$. This contradicts the definition of $T^{+}$. Hence $T^{+} = \infty $. 
A similar agrument shows that ${\mathcal{T}}^{-}$, the set of all $t \in {\R }_{<0}$ such that 
there is a continuous curve $c: (-t, 0] \subseteq \R \rightarrow M$, which covers 
the given continuous curve $\overline{c}$ on $(-t, 0 ]$, is nonempty and does not have 
a finite infimum. Consequently, for a given continuous curve $\overline{c}: \R \rightarrow 
\overline{M}$ there is a continuous curve $c: \R \rightarrow M$, which covers $\overline{c}$. 
\hfill $\square $ \medskip 

In our context we have 
\begin{displaymath}
{\Gamma }_0C^{\infty}(M)^G = \{ c \in C^0(\R , M) \setrule \, 
f \comp c \in C^{\infty}(\R )\, \, \mbox{for all $f \in C^{\infty}(M)^G$} \} 
\end{displaymath}
and 
\begin{displaymath}
{\Phi }_0{\Gamma }_0C^{\infty}(\overline{M}) = \{ \overline{f} \in C^0(\overline{M}) \setrule \, 
\overline{f} \comp \overline{c} \in C^{\infty}(\R ) \, \, 
\mbox{for all $\overline{c} \in {\Gamma }_0C^{\infty}(\overline{M})$} \} ,
\end{displaymath} 
see equation (\ref{eq-s3ss3threevnw}). \medskip 

\noindent \textbf{Lemma 6.2.4} We have 
\begin{equation}
{\pi  }_{\ast }({\Gamma }_0 C^{\infty}(M)^G)  = {\Gamma }_0C^{\infty}(\overline{M}),  
\label{eq-s5onenww}
\end{equation}
that is, if $\overline{c} \in {\Gamma }_0C^{\infty}(\overline{M})$, then there is a 
$c \in {\Gamma }_0C^{\infty}(M)^G$ such that $\overline{c} = {\pi }_{\ast }(c) = \pi \comp c$ and 
conversely, if $c \in {\Gamma }_0C^{\infty}(M)^G$, then ${\pi }_{\ast }(c) = \overline{c} 
\in {\Gamma }_0C^{\infty}(\overline{M})$. \medskip 

\noindent \textbf{Proof.} Suppose that $\overline{c} \in {\Gamma }_0C^{\infty}(\overline{M}) 
\subseteq C^0(\overline{M})$. Then $\overline{f} \comp \overline{c} \in C^{\infty}(\R )$ 
for every $\overline{f} \in C^{\infty}(\overline{M})$. By lemma 6.2.3 there is $c \in C^0(\R , M)$ such 
that $\overline{c} = {\pi }_{\ast }(c)$. Hence 
\begin{displaymath}
C^{\infty}(\R ) \ni \overline{f} \comp \overline{c} = (\overline{f} \comp \pi ) \comp c 
= ({\pi }^{\ast }\overline{f}) \comp c = f \comp c
\end{displaymath}
for every $f \in C^{\infty}(M)^G$, since $C^{\infty}(M)^G = {\pi }^{\ast }(C^{\infty}(\overline{M}))$. 
In other words, $c \in {\Gamma }_0C^{\infty}(M)^G$. So $\overline{c} = {\pi }_{\ast }(c) 
\in {\pi  }_{\ast }({\Gamma }_0C^{\infty}(M)^G)$. Consequently, ${\Gamma }_0C^{\infty}(\overline{M}) 
\subseteq {\pi }_{\ast }( {\Gamma }_0C^{\infty}(M)^G)$. To verify the reverse inclusion suppose 
$\overline{c} \in {\pi }_{\ast }( {\Gamma }_0C^{\infty}(M)^G)$. Then there is a 
$c \in {\Gamma }_0C^{\infty}(M)^G$ such that $\overline{c} = {\pi }_{\ast }(c)$. Suppose that 
$\overline{f} \in C^{\infty}(\overline{M})$. Then $f = {\pi }^{\ast }\overline{f} \in C^{\infty}(M)^G$. Now
\begin{displaymath}
C^{\infty}(\R ) \ni f \comp c = ({\pi }^{\ast }\overline{f}) \comp c = \overline{f}\comp (\pi \comp c) 
= \overline{f} \comp \overline{c}, 
\end{displaymath}
for every $f \in C^{\infty}(M)^G$ and hence for all $\overline{f} \in C^{\infty}(\overline{M})$, 
since $c \in {\Gamma }_0C^{\infty}(M)^G$ and ${\pi }^{\ast }(C^{\infty}(\overline{M})) = 
C^{\infty}(M)^G$. Thus $\overline{c} \in {\Gamma }_0C^{\infty}(\overline{M})$. So we have 
shown that ${\Gamma }_0C^{\infty}(M)^G \subseteq {\Gamma }_0C^{\infty}(\overline{M})$. 
This verifies (\ref{eq-s5onenww}). \hfill $\square $ \medskip 

\noindent \textbf{Claim 6.2.5} The differential structure $C^{\infty}(\overline{M})$ is 
continuously reflexive, that is 
\begin{equation}
C^{\infty}(\overline{M}) = {\Phi }_0{\Gamma }_0(C^{\infty}(\overline{M})). 
\label{eq-s5threenww}
\end{equation}

To prove claim 6.2.5 we will need some additional results. \medskip 

\noindent \textbf{Lemma 6.2.6} We have 
\begin{equation}
{\Gamma }_0C^{\infty}(M) \subseteq {\Gamma }_0C^{\infty}(M)^G.
\label{eq-s5threeAnw}
\end{equation}

\noindent \textbf{Proof.} Let $c \in {\Gamma }_0C^{\infty}(M)$. Then $c \in C^0(\R , M)$ and 
for every $f \in C^{\infty}(M)$ we have $f \comp c \in C^{\infty}(\R )$. Thus 
$f \comp c \in C^{\infty}(\R )$ for every $f \in C^{\infty}(M)^G$ because $C^{\infty}(M)^G 
\subseteq C^{\infty}(M)$. So $c \in {\Gamma }_0C^{\infty}(M)^G$. This verifies 
(\ref{eq-s5threeAnw}). \hfill $\square $ \medskip 

\noindent \textbf{Lemma 6.2.7} We have 
\begin{equation}
{\Phi }_0{\Gamma }_0C^{\infty}(M)^G \subseteq 
{\Phi }_0 {\Gamma }_0C^{\infty}(M) . 
\label{eq-s5threeBnw}
\end{equation}

\noindent \textbf{Proof.} Let $f \in {\Phi }_0{\Gamma }_0C^{\infty}(M)^G \subseteq 
C^0(M)^G$. Then $f \comp c \in C^{\infty}(\R )$ for all $c \in {\Gamma }_0C^{\infty}(M)^G$. 
Since ${\Gamma }_0C^{\infty}(M) \subseteq {\Gamma }_0C^{\infty}(M)^G$ by lemma 6.2.6, 
we obtain $f \comp c \in C^{\infty}(\R )$ for all $c \in {\Gamma }_0C^{\infty}(M)$. In 
other words, $f \in {\Phi }_0{\Gamma }_0C^{\infty}(M)$. This verifies  
(\ref{eq-s5threeBnw}). \hfill $\square $ \medskip 

\noindent \textbf{Proposition 6.2.8} $C^{\infty}(M)^G$ is continuously reflexive. \medskip %

\noindent \textbf{Proof.} Using the fact that $C^{\infty}(M)^G \subseteq 
{\Phi}_0{\Gamma }_0(C^{\infty}(M)^G)$, we need only show that 
${\Phi}_0{\Gamma }_0(C^{\infty}(M)^G) \subseteq C^{\infty}(M)^G$. Suppose that 
$f \in {\Phi }_0{\Gamma }_0C^{\infty}(M)^G $. Then $f \in C^0(M)^G$, that is, 
$f$ is a continuous $G$-invariant function on $M$. Since ${\Phi }_0{\Gamma }_0C(M)^G 
\subseteq {\Phi }_0{\Gamma }_0C^{\infty}(M)$ by lemma 6.2.7 and 
$f \in {\Phi }_0{\Gamma }_0C^{\infty}(M)^G$ by \linebreak 
hypothesis, we get 
$f \in {\Phi }_0{\Gamma }_0C^{\infty}(M) = C^{\infty}(M)$. The preceding equality 
follows since the differential structure $C^{\infty}(M)$ on the smooth manifold $M$ is continuously reflexive. So $f$ is a smooth function on $M$. Since $f$ is $G$-invariant, it follows that 
$f \in C^{\infty}(M)^G$. This shows that ${\Phi }_0{\Gamma }_0C^{\infty}(M)^G \subseteq 
C^{\infty}(M)^G$. \hfill $\square $ \medskip 

\noindent \textbf{Lemma 6.2.9} We have 
\begin{equation}
{\pi }^{\ast } \big( {\Phi }_0{\Gamma }_0C^{\infty}(\overline{M}) \big) 
= {\Phi }_0{\Gamma }_0C^{\infty}(M)^G .
\label{eq-s5threeCnw}
\end{equation}

\noindent \textbf{Proof.} Suppose that $\overline{f} \in {\Phi }_0{\Gamma }_0C^{\infty}(\overline{M}) 
\subseteq C^0(\overline{M})$. Then for every $\overline{c} \in {\Gamma }_0C^{\infty}(\overline{M})$ 
we have $\overline{f} \comp \overline{c} \in C^{\infty}(\R )$. Since 
${\Gamma }_0C^{\infty}(\overline{M}) = {\pi }_{\ast }( {\Gamma }_0C^{\infty}(M)^G)$ by lemma 
6.2.4, there is $c \in {\Gamma }_0C^{\infty}(M)^G$ such that $\overline{c} = \pi \comp c$. So 
\begin{displaymath}
C^{\infty}(\R ) \ni \overline{f} \comp \overline{c} = (\overline{f} \comp \pi ) \comp c 
={\pi }^{\ast }\overline{f} \comp c, 
\end{displaymath}
for every $c \in {\Gamma }_0C^{\infty}(M)^G$. Thus ${\pi }^{\ast }\overline{f} 
\in {\Phi }_0{\Gamma }_0C^{\infty}(M)^G$. So we have shown 
${\pi }^{\ast }\big( {\Phi }_0{\Gamma }_0C^{\infty}(\overline{M}) \big) \subseteq 
{\Phi }_0{\Gamma }_0C^{\infty}(M)^G$. To verify the reverse inclusion let 
$f \in {\Phi }_0{\Gamma }_0C^{\infty}(M)^G \subseteq C^0(M)^G$. Then 
$f \comp c \in C^{\infty}(\R )$ for every $c \in {\Gamma }_0C^{\infty}(M)^G$. Since 
$C^0(M)^G = {\pi }^{\ast }(C^0(\overline{M}))$, there is $\overline{f} \in C^0(\overline{M})$ 
such that $f = {\pi }^{\ast }\overline{f}$. So 
\begin{displaymath}
C^{\infty}(\R ) \ni f \comp c = {\pi }^{\ast }\overline{f} \comp c = \overline{f} \comp (\pi \comp c), 
\end{displaymath}
for every $c \in {\Gamma }_0C^{\infty}(M)^G$. By lemma 6.2.4, 
${\pi  }_{\ast } ({\Gamma }_0C^{\infty}(M)^G) = {\Gamma }_0C^{\infty}(\overline{M})$. Consequently, $\overline{f} \comp (\pi \comp c) = 
\overline{f} \comp \overline{c} \in C^{\infty}(\R )$ for every $\overline{c} \in 
{\Gamma }_0C^{\infty}(\overline{M})$. So $\overline{f} \in 
{\Phi }_0{\Gamma }_0C^{\infty}(\overline{M})$. In other words, $f \in {\pi }^{\ast } \big( 
{\Phi }_0{\Gamma }_0C^{\infty}(\overline{M}) \big) $. Thus  
$ {\Phi }_0{\Gamma }_0C^{\infty}(M)^G  \subseteq 
{\pi }^{\ast} \big( {\Phi }_0{\Gamma }_0C^{\infty}(\overline{M}) \big) $. This verifies 
(\ref{eq-s5threeCnw}). \hfill $\square $ \medskip 

\noindent \textbf{Proof of claim 6.2.5} From lemma 6.2.9 it follows that we have 
${\pi }^{\ast }\big( {\Phi }_0{\Gamma }_0C^{\infty}(\overline{M})\big) $ 
$= {\Phi }_0{\Gamma }_0C^{\infty}(M)^G$. But $C^{\infty}(M)^G$ is continuously reflexive 
by proposition 6.2.8. Thus 
\begin{displaymath} 
{\pi }^{\ast }\big( {\Phi }_0{\Gamma }_0C^{\infty}(\overline{M})\big) = C^{\infty}(M)^G = 
{\pi }^{\ast }(C^{\infty}(\overline{M})) . 
\end{displaymath}
Since the $G$-orbit mapping $\pi $ is surjective, the map 
\begin{displaymath}
{\pi }^{\ast }: C^{\infty}(\overline{M}) \rightarrow C^{\infty}(M)^G: \overline{f} \mapsto 
{\pi }^{\ast }(\overline{f})
\end{displaymath}
is injective. Hence ${\Phi }_0{\Gamma }_0C^{\infty}(\overline{M}) = C^{\infty}(\overline{M})$, 
which verifies (\ref{eq-s5threenww}). \hfill $\square $ \medskip  

\noindent \textbf{Proof of theorem 6.2.2} This follows from theorem 4.3, since the 
differential structure $C^{\infty}(\overline{M})$ is continuously reflexive by claim 6.2.5. 
\hfill $\square $ 

\subsection{Additional properties}

We now verify the extra hypotheses on a differentiable space  
so that the \linebreak 
hypotheses of Smith's de Rham theorem, theorem 5.5.3, hold for 
the differentiable space $(\overline{M}, C^{\infty}(\overline{M}))$. \medskip

\noindent \textbf{Claim 6.3.1} The subcartesian differential space 
$(\overline{M}, C^{\infty}(\overline{M}))$ is a locally differentiably contractible 
differentiable space. \medskip 

\noindent \textbf{Proof.} Let $m \in M$. Since the $G$-action $\Phi $ on the smooth 
manifold $M$ is proper, there is a $G_m$-invariant slice $S$ to the $G$-action 
at $m$. Because the isotropy group $G_m = \{ g \in G \setrule \, {\Phi }_g(m) = m \}$ is compact, there is an open neighborhood $V_m$ of $0$ in $T_mM$ and a diffeomorphism 
$\psi : V_m \subseteq T_mM \rightarrow 
U_m \subseteq M$ such that for every $g \in G_m$ and every $v_m \in V_m$ we have 
$T_m\psi (T_m {\Phi }_g v_m) = T_m {\Phi }_g (T_m \psi \, v_m)$. Let ${\gamma }_m$ be a 
$G_m$-invariant Euclidean inner product on $T_m M$ and let $B^{\varepsilon }_m$ be 
a ${\gamma }_m$ ball of radius $\varepsilon $ with center $0$ contained in $V_m$. Let 
$Y_m$ be a linear vector field on $T_mM$ whose flow is ${\varphi }^{Y_m}_t(v_m) = 
{\mathrm{e}}^{-t}v_m$. Let $X_m = \frac{1}{\mathrm{vol}G_m} \int_{G_m} 
(T_m{\Phi }_{h})^{\ast }Y_m \, \dee h$, where $\dee h$ is a volume form on $G_m$, be the vector 
field on $T_mM$, which is the $G_m$ average of the pull back of $Y_m$ by the action 
$T_m{\Phi }_h$ of the $G_m$ action on $T_m$. Then $X_m$ is a 
linear $G_m$-invariant vector field on $T_mM$, whose flow ${\varphi }^{X_m}_t$ commutes 
with the ${\varphi }^{Y_m}_t$ of the vector field $Y_m$. Hence the flow of the $G_m$-invariant vector field ${\mathcal{X}}_m = {\psi }_{\ast }X_m$ on $U_m$ preserves the slice 
$\mathcal{S} = S \cap \psi (B^{\varepsilon }_m)$ and contracts it to the point $m$. 
The vector field $\mathcal{X}({\Phi }_g(m)) = ({\Phi }_g)_{\ast }{\mathcal{X}}_m$ is 
${\Phi }_g$-invariant on $G \cdot \mathcal{S}$. Its flow ${\varphi }^{\mathcal{X}}_t$ 
contracts $G\cdot \mathcal{S}$ to the $G$-orbit $G_m$ and commutes with the $G$ 
action $\Phi $. Let $\overline{\mathcal{S}} = \pi (G \cdot \mathcal{S})$. Then 
$\overline{\mathcal{S}}$ is an open neighborhood of $\overline{m} = \pi (m)$ in the 
quotient topology (and hence the differential space topology) of 
$(\overline{M}, C^{\infty}(\overline{M}))$. Let $\overline{\mathcal{X}}$ be the vector field on 
$\overline{M}$ induced by the $G$-invariant vector $\mathcal{X}$. Since 
${\varphi }^{\overline{\mathcal{X}}} \comp \pi  = \pi \comp 
{\varphi }^{\mathcal{X}}_t$, the flow ${\varphi }^{\overline{\mathcal{X}}}$ is $\pi $-related 
to the flow ${\varphi }^{\mathcal{X}}_t$ of $\mathcal{X}$ on $M$, that is, 
${\varphi }^{\overline{\mathcal{X}}}_t \comp \pi = \pi \comp {\varphi }^{\mathcal{X}}_t$. 
From the fact that ${\varphi }^{\mathcal{X}}_t$ smoothly contracts $G\cdot \mathcal{S}$ to the 
$G$-orbit $G\cdot m$, we determine that ${\varphi }^{\overline{\mathcal{X}}}_t$ smoothly 
contracts the open neighborhood $\overline{\mathcal{S}}$ of $\overline{m}$ to 
$\overline{m}$. \hfill $\square $ \medskip 

Let $(S, C^{\infty}(S))$ be a locally closed (and hence paracompact) subcartesian space. Let $U$ be an open subset of $S$. The assignment 
\begin{displaymath}
U \mapsto C^{\infty}(U) = C^{\infty}(S)|_U
\end{displaymath}
defines an abelian sheaf $\mathcal{S}$ of rings with a unit. \medskip 

\noindent \textbf{Lemma 6.3.2} The sheaf $\mathcal{S}$ is \emph{soft}, that is, 
if $T$ is a closed subset of $S$ and $s \in C^{\infty}(T) = C^{\infty}(S)|_T$, then 
there is a $\Sigma \in C^{\infty}(S)$ such that $\Sigma |_T = s$. \medskip 

\noindent \textbf{Proof.} Let $U_0 = S \setminus T$ and let ${\{ U_i \}}_{i \in I}$ be an 
open covering of $S$. Then ${\{ U_i \} }_{j \in I\cup \{0 \}}$ is an open covering of 
$S$. Since $S$ is paracompact, ${\{ U_i \} }_{j \in I\cup \{0 \}}$ has a locally 
finite subcovering $\mathcal{V} = {\{ U_j \} }_{j \in J \subset (I \cup \{ 0 \}) }$, 
Let ${\{ {\psi }_j \} }_{j\in J}$ be a partition of unity in $C^{\infty}(S)$ associated to 
$\mathcal{V}$. Then $\Sigma = \sum_{j \in J} {\psi }_j s$ is a smooth function on $S$, 
whose restriction to $T$ is $s$, because $\sum_{j \in J} {\psi }_j(p) =1$ for every 
$p \in T$. Hence the sheaf $\mathcal{S}$ is soft. \hfill $\square $ \medskip 
 
\noindent \textbf{Theorem 6.3.3} The paracompact subcartesian space 
$(\overline{M}, C^{\infty}(\overline{M}))$, where $\overline{M}$ is the orbit space 
of a proper action of a connected Lie group on a smooth manifold 
$M$ is a locally closed differential space, which is locally differentiably contractible and 
whose sheaf $\mathcal{S}$ of smooth functions is soft. \medskip 

\noindent \textbf{Proof.} This follows immediately from claim 6.3.1 and lemma 6.3.2. 
\hfill $\square $ 

\subsection{de Rham's theorem for $(\overline{M}, C^{\infty}(\overline{M})$}

Using theorem 6.3.3 it follows that Smith's version of de Rham's theorem 5.5.3 holds for 
$(\overline{M}, C^{\infty}(\overline{M}))$, namely, \medskip 

\noindent \textbf{Theorem 6.4.1} The cohomology of the exterior differential algebra 
$(\Lambda (\overline{M}), \overline{\dee}\, )$ of exterior differential forms on 
$(\overline{M}, C^{\infty}(\overline{M}))$ is isomorphic to the 
subalgebra $\mathcal{A}$ of the exterior differential algebra $( {\Lambda (M)}^G,\dee \, )$ of 
$G$-invariant differential forms on $M$ generated by differential $1$-forms of 
smooth $G$-invariant functions on $M$ with coefficients which are 
smooth $G$-invariant functions. \medskip 

The following example shows that the differential forms in Smith's version of 
de Rham's theorem \cite{smith} are not the basic forms used in Koszul's formulation 
of de Rham's theorem, see \cite{koszul} and also Karshon and Watts \cite{karshon-watts}. \medskip 

\noindent \textbf{Example 6.4.2} We follow Smith \cite[p.133--134]{smith}. Consider the 
${\Z }_2 = \{ -1, 1 \} $ action on ${\R }^2$ with coordinates $(x,y)$ generated by 
\begin{equation}
\zeta : {\R }^2 \rightarrow {\R }^2: (x,y) \mapsto (-x,-y). 
\label{eq-s5one}
\end{equation}
Since ${\Z }_2$ is finite, the ${\Z}_2$ action on ${\R }^2$ is proper. It is straightforward 
to check that the algebra of ${\Z }_2$ invariant polynomials on ${\R }^2$ is 
generated by 
\begin{equation}
{\pi }_1 = x^2, \, \, \, {\pi }_2 = xy, \, \, \mathrm{and}\, \, {\pi }_3 = y^2. 
\label{eq-s5two}
\end{equation}
They are subject to the relation 
\begin{equation}
{\pi }^2_3 = {\pi }_1{\pi }_2, \, \, {\pi }_1 \ge 0 \, \, \& \, \, {\pi }_2 \ge 0, 
\label{eq-s5twoA}
\end{equation}
which defines the orbit space ${\R }^2/{\Z}_2$. \medskip 

Consider the differential $1$-form $\omega = x \dee y - y \dee x$ on ${\R }^2$, 
which is clearly ${\Z }_2$ invariant. Since the Lie algebra {\footnotesize ${\Z }_2$} of 
${\Z }_2$ is $\{ 0 \} $, we have $L_{X_{\xi }}\omega =0$ for every 
$\xi \in \mbox{\footnotesize ${\Z }_2$} $. So $\omega $ is a \emph{basic} $1$-form on 
${\R }^2$ and hence induces a $1$-form $\overline{\omega }$ on ${\R }^2/{\Z}_2$. \medskip 

We now show that $\omega $ is not a sum of $1$-forms of the form 
$f_0 \dee f_1$, where $f_0$ is a smooth ${\Z }_2$ invariant function on ${\R }^2$ and 
$f_1$ is a ${\Z }_2$ invariant polynomial. Using polar coordinates $( \theta , r)$ on ${\R }^2$, we see that $\omega = r^2 \dee \theta $. 
Let ${\gamma }_R$ be the circle $x^2+y^2 = R^2$ in ${\R }^2$. Then 
\begin{equation}
\int_{{\gamma }_R} \omega = 2\pi \, R^2 = \mathrm{O}(R^2).
\label{eq-s5three}
\end{equation}
Since the invariant polynomials ${\pi }_j$ are quadratic in $x$ and $y$, all the 
odd order terms in the Taylor expansion of $f_i$ at $(0,0)$ to fourth order vanish. 
Hence 
\begin{displaymath}
f_i(x,y) = a_i + b_i x^2 +c_i (xy) + d_iy^2 + \mathrm{O}((x^2+y^2)^2) ,  
\end{displaymath}
where $a_i$, $b_i$, $c_i$, and $d_i \in \R $. Thus 
\begin{displaymath}
\int_{{\gamma }_R} f_0 \dee f_1 = 
b_0c_1 \int_{{\gamma }_R} x^2 \dee \, (xy) + d_0c_1 \int_{{\gamma }_R} y^2 \dee \, (xy) 
+ \mathrm{O}\big( \int_{{\gamma }_R} r^4 \big) ,
\end{displaymath}
since a straightforward calculation shows that the integral $\int_{{\gamma }_R}$ of 
the $1$ forms $x\dee x$, $\dee \, (xy)$, $y\dee y$, $x^3\dee x$, $x^2y \dee y$, 
$(xy)\dee \, (xy)$, $y^2\dee x$, and $y^3 \dee y$ vanish. But 
\begin{displaymath}
\int_{{\gamma }_R} x^2 \dee \, (xy) = -R^4 \int^{2\pi }_0 {\cos }^2\theta \, \dee \theta 
+ R^4  \int^{2\pi }_0 {\cos }^4\theta \, \dee \theta = CR^4,
\end{displaymath}
where $C \ne 0$ and $\int_{{\gamma }_R} y^2 \dee \, (xy) =
 -\int_{{\gamma }_R} x^2 \dee \, (xy) $. 
Hence 
\begin{equation}
\int_{{\gamma }_R} f_0 \dee f_1 = c_1C(b_0-d_0)R^4 + \mathrm{O}(R^5) = 
\mathrm{O}(R^4), 
\label{eq-s5four}
\end{equation}
which contradicts equation (\ref{eq-s5three}). \medskip %

Thus the exterior algebra of basic forms on ${\R }^2$ differs from the exterior 
algebra generated by the exterior derivative of ${\Z }_2$ invariant functions 
with coefficients which are smooth ${\Z }_2$-invariant functions. However, Koszul's 
version of de Rham's theorem and Smith's version yield the same real cohomology groups for 
${\R }^2/{\Z }_2$, because they both give \v{C}ech cohomology with values in the sheaf $R$ of locally 
constant real valued functions on ${\R }^2/{\Z }_2$, since ${\R }^2$ is connected. \hfill $\square $

\end{document}